\documentclass[tbtags]{m2anJWB}

\usepackage{verbatim}
\usepackage{amssymb,amsmath,psfrag}
\usepackage{pstricks,ifthen}
\usepackage{color}
\usepackage{bm}
\usepackage{floatrow}

\usepackage{graphicx}
\usepackage{multirow}
\allowdisplaybreaks

\def\epsilon{\varepsilon}
\def\hat{\widehat}

\newfloatcommand{capbtabbox}{table}[][\FBwidth]
\newcommand{\mintpis}{\int_{\partial_i \sigma}\!\!\!\!\!\!\!\!\!\!\!{\rm-}\ }
\newcommand{\mints}{\int_{\sigma}\!\!\!\!\!\!\!{\rm-}\ }
\newcommand{\mintD}{\int_{\,D}\!\!\!\!\!\!\!\!\!{\rm-}\ \ }

\newcommand{\va}{\underline{a}}

\newcommand{\vd}{\underline{d}}
\newcommand{\ve}{\underline{e}}

\newcommand{\vj}{\underline{j}}
\newcommand{\vJ}{\underline{J}}

\newcommand{\vpsi}{\underline{\psi}}

\newcommand{\vq}{\underline{q}}
\newcommand{\vQ}{\underline{Q}}
\newcommand{\vB}{\underline{B}}
\newcommand{\vP}{\underline{P}}

\newcommand{\vI}{\underline{I}}

\newcommand{\vv}{\underline{v}}

\newcommand{\vx}{\underline{x}}
\newcommand{\vy}{\underline{y}}

\newcommand{\vz}{\underline{z}}
\newcommand{\vZ}{\underline{Z}}

\newcommand{\vzero}{\underline{0}}
\newcommand{\vV}{\underline{V}}

\newcommand{\vVh}{\underline{V}^h}

\newcommand{\Hhalf}{H^{\frac{1}{2}}_{00}(\Omega)}
\newcommand{\del}{\underline{\nabla}}
\newcommand{\eps}{\varepsilon}
\newcommand{\ts}{\textstyle}
\newcommand{\ds}{\displaystyle}
\newcommand{\dx}{\,{\rm d}\vx}
\newcommand{\dy}{\,{\rm d}\vy}
\newcommand{\dt}{{\rm d}t}

\newcommand{\beq}{\begin{equation}}
\newcommand{\eeq}{\end{equation}}

\renewcommand{\theequation}{\arabic{section}.\arabic{equation}}

\newcounter{ind}
\def\eqlabstart{%
\setcounter{ind}{\value{equation}}\addtocounter{ind}{1}%
\setcounter{equation}{0}%
\renewcommand{\theequation}{\arabic{section}.\arabic{ind}\alph{equation}}}
\def\eqlabend{%
\renewcommand{\theequation}{\arabic{section}.\arabic{equation}}%
\ifthenelse{\value{equation}=0}{\addtocounter{ind}{-1}}{}%
\setcounter{equation}{\value{ind}}}

\newcommand{\veta}{\underline{\eta}}
\newcommand{\vc}{\underline{c}}
\newcommand{\vS}{\underline{S}}
\newcommand{\vG}{\underline{G}}

\begin{document}
\title{ Sandpiles and Superconductors:\\ Nonconforming Linear Finite Element
Approximations for Mixed Formulations of Quasi-Variational Inequalities
}
\author{John W. Barrett}
\address{Department of Mathematics,
Imperial College London, London, SW7 2AZ, UK.}

\author{Leonid Prigozhin}
\address{Department of Solar Energy and Environmental Physics,
Blaustein Institutes for Desert Research,
Ben-Gurion University of
the Negev, Sede Boqer Campus, 84990 Israel.}

\runningauthors{John W. Barrett and Leonid Prigozhin}

\runningtitle{Quasi-Variational Inequality Problems Arising in the
Modeling of Growing Sandpiles}

\date{}

\begin{abstract}
Similar evolutionary variational and quasi-variational inequalities
with gradient constraints arise in the modeling of growing sandpiles
and type-II superconductors.
Recently, mixed formulations of these inequalities were used for establishing existence results
in the quasi-variational inequality case. Such formulations, and this is an additional advantage,
made it possible to determine numerically not only the primal variables, e.g.\ the evolving sand
surface and the magnetic field for sandpiles and superconductors, respectively,
but also the dual variables, the sand flux and the electric field.

Numerical approximations of these mixed formulations in previous works employed
the Raviart--Thomas element of the lowest order.
Here we introduce simpler numerical approximations of these mixed formulations
based on the nonconforming linear finite element. We prove (subsequence) convergence
of these approximations,
and illustrate their effectiveness by numerical experiments.
\end{abstract}


\vspace{2mm}


\maketitle

\section{Introduction}
\label{intro}
\setcounter{equation}{0}

Recently, the present authors have introduced mixed formulations
of variational and quasi-variational inequality problems
arising in the mathematical modelling of (i) growing sandpiles, (ii) cylindrical superconductors
in a parallel external field and (iii) thin film superconductors in a perpendicular
external field in \cite{sandqvi}, \cite{BP5} and \cite{EFILMMATH},
respectively. In each of these papers, a numerical approximation,
based on the lowest order Raviart--Thomas element,
of the corresponding mixed formulation was introduced, and (subsequence) convergence was proved
as the mesh parameters and the power law regularization parameter, $r-1$,
tended to zero.
Hence, the existence of a solution to these mixed formulations was established.
In this paper, we introduce simpler numerical approximations
based on a nonconforming linear finite element approximation
of these mixed formulations. In addition, we
prove (subsequence) convergence of these approximations
as the mesh and regularization parameters tend to zero.

We first briefly describe these mixed formulations.
Let $\Omega \subset \mathbb{R}^2$ be a simply connected domain
with a Lipschitz boundary $\partial \Omega$.

\subsection{Mathematical models and their mixed formulations}

\subsubsection*{(i) Growing Sandpiles}
Let a cohesionless granular material (sand), characterized by its angle of repose $\alpha$,
be poured out onto a rigid surface $y=w_0(\vx)$, where $y$ is vertical and $\vx\in\Omega$.
The support surface
$w_0\in W^{1,\infty}_0(\Omega)$
and the nonnegative density of the distributed source  $f \in L^2(0,T;L^2(\Omega))$ are given.
We consider the growing sandpile $y=w(\vx,t)$
and set an open boundary condition $w|_{\partial\Omega}=0$.
Denoting by $\vq(\vx,t)$ the horizontal projection of the flux of material pouring
down the evolving pile surface,
we can write the mass balance equation
\beq \frac{\partial w}{\partial t} +\del \,.\, \vq=f.\label{bal}\eeq
The quasi-stationary model of sand surface evolution, see Prigozhin \cite{Fill,River,Psand},
assumes the flow
of sand is confined to a thin surface layer and directed towards the steepest descent of
the pile surface. Wherever the support surface is covered by sand, the pile slope should not
exceed the critical value; that is, $w>w_0\ \Rightarrow\ |\del w|\le k_0$,
where $k_0=\tan\alpha$
is the internal friction coefficient. Of course,
the uncovered parts of the support can be steeper.
This model does not allow for any flow on the subcritical parts of the pile surface; that is,
$|\del w|< k_0\ \Rightarrow\ \vq=\underline{0}$.
These constitutive relations can be conveniently reformulated
for a.e.\ $(\vx,t) \in \Omega \times (0,T)$ as
\begin{align}|\del w| \leq M(w)
\qquad \mbox{and} \qquad
M(w)\,|\vq| +\del w\,.\,\vq =0,\label{cond1}
\end{align}
where, for any $\eta\in C(\overline{\Omega})$,
\begin{align}
M(\eta)(\vx)&:=\left\{ \begin{array}{ll}  k_0 \qquad &\eta(\vx)>w_0(\vx),\\
\max(k_0,|\del w_0(\vx)|) \qquad & \eta(\vx)\le w_0(\vx)\end{array}\right.
\qquad \forall \ \vx \in \overline{\Omega}.
\label{M}
\end{align}
Let us define, for any $\eta\in C(\overline{\Omega})$, the closed convex non-empty set
\beq K(\eta):=\{ \varphi\in W^{1,\infty}_0(\Omega)\ :\ |\del\varphi|\le M(\eta)\ \
\mbox{a.e. in}\ \Omega\}. \label{Kset}\eeq
Since $M(w)\,|\vq|+\del \varphi \,.\, \vq \geq 0$ for any $\varphi\in K(w)$,
we have, on noting (\ref{cond1}), that $w \in K(w)$ and $\del(\varphi-w)\,.\,\vq\geq 0$.
A weak form of the latter inequality is:
for a.a. $t \in (0,T)$
\beq
\int_{\Omega} \del\,.\,\vq\,(w-\varphi)\, {\rm d}\vx \geq 0
\qquad \forall \ \varphi \in K(w).
\label{sbd1}
\eeq
Combining (\ref{sbd1}) and (\ref{bal}) yields an evolutionary quasi-variational
inequality for the evolving pile surface: Find $w \in K(w)$ such that
for $a.a.\ t \in (0,T)$
\beq
\int_{\Omega} \left(\frac{\partial w}{\partial t} - f \right)(\varphi-w)\,
{\rm d}\vx \geq 0
\qquad \forall \ \varphi \in K(w).
\label{P1}
\eeq
Assuming there is no sand  on the support initially,  we set
\beq w(\cdot,0)=w_0(\cdot). 
\label{ic}\eeq
We note that 
with the open boundary condition
$w|_{\partial\Omega}=0$ an uncontrollable influx of material from outside can occur
through the parts of the boundary where $\del w\,.\,\underline{\nu}\geq k_0$,
with $\underline{\nu}$ being the outward unit normal to $\partial \Omega$.
This makes the solution non-unique and, possibly, discontinuous.
Such an influx is prevented in our model by assuming that
\beq\del w_0\,.\,\underline{\nu}< k_0\ \ \mbox{on}\
\partial\Omega,\label{influx}\eeq 
which implies, see \cite{sandqvi}, that $\del w\,.\,\underline{\nu}< k_0$ on
$\partial\Omega$ also for $t>0$.

If $|\del w_0| \leq k_0$ a.e. in $\Omega$, then
$K(\eta)\equiv K:=\{ \varphi\in W^{1,\infty}_0(\Omega)\ :\ $
$|\del\varphi|\le k_0\ \ \mbox{a.e. in}\ \Omega\}$ and the quasi-variational inequality
(\ref{P1}) becomes simply a variational inequality; this case was studied in
Prigozhin \cite{Fill,Psand} and Aronson, Evans and Wu \cite{AEW}.

Here we will use a mixed variational formulation 
of the growing  sandpile model involving both variables.
Such formulations are often advantageous,
because they allow one to determine not only the evolving sand surface $w$ but also the
surface flux $\vq$, which is of interest too in various applications; see Prigozhin \cite{PChEJ,River},
and Barrett and Prigozhin \cite{BP5}.
In such formulations, and this is their additional advantage, the difficult to deal with
gradient constraint in (\ref{P1}) is replaced by a simpler, although non-smooth, nonlinearity.
Therefore instead of excluding the surface flux $\vq$ from
the model formulation, as in the transition to (\ref{P1}) above,
we 
 reformulate the conditions 
(\ref{cond1}) for a.a.\ $t \in (0,T)$ as
\beq \int_{\Omega}\left[M(w)\,(|\vv|-|\vq|)+ \del w\,.\,(\vv-\vq)\right] {\rm d}\vx\geq 0\label{cond3}\eeq
for any test flux $\vv$, and consider a mixed formulation of the sand model
as (\ref{bal}) and (\ref{cond3}).

The natural function space for the flux $\vq$ is the
space of vector-valued bounded Radon measures having $L^2$ divergence. If $\vq$ is such a
measure,  the discontinuity of $M(w)$ makes it  difficult to give a sense to
the term $\int_{\Omega}M(w)\, |\vq| \,{\rm d}\vx $
in the inequality (\ref{cond3}) of the mixed formulation.
Existence of a solution was recently proved in Barrett and Prigozhin \cite{sandqvi},
for a regularized version of the growing sandpile model
with a
continuous operator $M_\epsilon : C(\overline{\Omega}) \rightarrow C(\overline{\Omega})$,
determined as follows.
For a fixed small $\varepsilon >0$, we approximate the initial data
$w_0\in W^{1,\infty}_0(\Omega)$ by $w_0^\epsilon \in W^{1,\infty}_0(\Omega) \cap
C^1(\overline{\Omega})$, and
$M(\cdot)$ by the continuous function $M_\eps(\cdot)$ such that
for any $\vx \in \overline{\Omega}$
\begin{align}
&M_{\varepsilon}(\eta)(\vx)
:=\left\{ \begin{array}{ll}  k_0 \quad& \eta(\vx)\geq w_0^\epsilon(\vx)
+ \varepsilon,\\
k_1^\epsilon(\vx) + (k_0-k_1^\epsilon(\vx)) \, \left(\ds\frac{\eta(\vx)-w_0^\epsilon(\vx)}{\varepsilon}\right)
\quad & \eta(\vx) \in [w_0^\epsilon(\vx),w_0^\epsilon(\vx)+\varepsilon],\\[4mm]
k_1^\epsilon(\vx)
:= \max(k_0,|\del w_0^\epsilon(\vx)|)
\quad & \eta(\vx)\le w_0^\epsilon(\vx).\end{array}\right.
\label{Meps}
\end{align}
Below, we also adopt such a regularisation. We note that the existence 
of a solution for the regularized
primal quasi-variational inequality (\ref{P1})
follows also from a recent result by Rodrigues and Santos \cite{RS}.

Obviously, if $|\nabla w^0|\leq k_0$ no regularisation is needed as $M\equiv k_0$.
In this variational inequality case the mixed formulations of the growing sandpile problem,
and its numerical approximation by the lowest order Raviart--Thomas element,
were studied in Barrett and Prigozhin \cite{Dual} and Dumont and Igbida \cite{DIa}.

\subsubsection*{(ii) Cylindrical Superconductors in a Parallel External Field}

Let us consider an infinite type-II superconducting cylinder having a cross section $\Omega$
and placed into a given parallel non-stationary uniform external magnetic field $b_e(t)$.
In this case the magnetic field of a current induced in the superconductor has also
only one non-zero component
and can be regarded as a scalar function $w(\vx,t)$, which vanishes on $\partial \Omega$.
The electric field, $\ve$,
inside the superconductor is the same in each cross section of the cylinder
and is orthogonal to the magnetic field.  A similar statement holds for the current density, $\vj$,
inside the superconductor.
With $\ve(\vx,t) \equiv [e_1(\vx,t),e_2(\vx,t)]^{\top}$,
Faraday's law can be rewritten as (\ref{bal}) with
\begin{align}
f= - \frac{d b_e}{dt} \qquad  \mbox{and} \qquad \vq = [e_2,-e_1]^{\top} \quad
\Rightarrow \quad
\nabla \times \ve = \frac{\partial e_2}{\partial x_1}-\frac{\partial e_1}{\partial x_2}
= \del\,.\,\vq\,.
\label{curldiv}
\end{align}
Here, and throughout this paper, we use scaled dimensionless electromagnetic variables.
In particular, we do not distinguish between the magnetic induction and the magnetic field
on assuming that the magnetic permeability of the superconductor is equal to that of a vacuum
and is scaled to unity.

Amp\`{e}re's law yields that the current density $\vj = \del \times w
= [\frac{\partial w}{\partial x_2},-\frac{\partial w}{\partial x_1}]^{\top}$,
and so $|\vj| = |\del w|$.
Let $\vj$ and $\ve$ satisfy the critical state model relations:
\beq
|\vj|\leq j_c,\qquad \quad |\vj|< j_c \quad \Rightarrow \quad \ve =\vzero,
\qquad \quad \ve \neq \vzero \quad \Rightarrow \quad
\ve \parallel \vj \quad \Rightarrow \quad
\vq \parallel -\del w\,,
\label{crit}
\eeq
where $j_c$ is the critical current density,
which may be constant or depend only on $\vx$
(the Bean model, see \cite{Bean})
or depend also on the total magnetic field, $w+b_e$  (the Kim model, see \cite{Kim}).
Similarly to the growing sandpile problem, one can show that $w$
satisfies the quasi-variational inequality problem (\ref{P1}) with $f$ as in (\ref{curldiv})
and $K(w)$ replaced by $\hat K(w+b_e)$, where
\begin{align}
\hat K(\psi) := \{ \eta \in W^{1,\infty}_0(\Omega)
: |\del \eta| \leq j_c(\psi)
\mbox{ a.e. in } \Omega
\}\,.
\label{hatK}
\end{align}
This is supplemented with $w(\cdot,0)=w_0(\cdot)$, where $w_0 \in \hat K(w_0+b_e(0))$.
Once again, if $j_c$ is independent of the total magnetic field, i.e.\ the Bean model,
this quasi-variational inequality problem collapses to a variational inequality problem.
Similarly, the conditions (\ref{crit}) can be reformulated as (\ref{cond3})
with $M(w)$ replaced by $j_c(w+b_e)$, and this supplemented with
(\ref{bal}) yields the mixed formulation of this cylindrical superconductor
problem, see Barrett and Prigozhin \cite{BP5} for further details.
We note that  $\vq=[-e_2,e_1]^\top$ in \cite{BP5}, see page 684 there.

In \cite{BP5}, and in this paper,
we assume for the critical state model that
\begin{align}
j_c(w+b_e)(\vx,t) = k(\vx)\,\hat M(w(\vx,t)+b_e(t)), \label{jcsc}
\end{align}
where
${\widehat M} \in
C({\mathbb R},[{\widehat M}_0,{\widehat M}_1])$
{with} $ {\widehat M}_0,\,{\widehat M}_1 \in {\mathbb R}_{>0},$
and $k \in C(\overline{\Omega})$ with $k(\vx) \geq k_{\rm min} >0$
for all $\vx \in \overline{\Omega}$.
In \cite{BP5} we exploited the fact that $|\del (w+b_e)| \leq k\, \widehat M(w+b_e)$
can be rewritten as
$|\del [\widehat F(w+b_e)]| \leq k$, where $\widehat F'(\cdot)=
[{\widehat M}(\cdot)]^{-1}$ and $\widehat F(0)=0$.
Clearly, such a reformulation is not applicable to $M(\cdot)$, (\ref{M}),
or $M_\epsilon(\cdot)$, (\ref{Meps}), for the growing sandpile problem.

Engineers often describe the current-voltage relation of type-II superconductors by a power law
\begin{align}
\ve = \left(\frac{|\vj|}{j_c}\right)^{p-2} \frac{\vj}{j_c}
\quad \Rightarrow \quad
\del w = -j_c\, |\vq|^{r-2} \vq\,, \qquad \mbox{where } \frac{1}{r} +\frac{1}{p}=1,
\label{power}
\end{align}
with the power $p$ typically between $10$ and $100$. 
As is well-known, the critical state model relations (\ref{crit})
can be regarded as the $p \rightarrow \infty$ ($r \rightarrow 1$) limit of
the power law (\ref{power}); see Barrett and Prigozhin \cite{NA} in the case
of the homogeneous Bean model, $j_c \in \mathbb{R}_{>0}$, and Theorem
\ref{conthmsc} below for (\ref{jcsc}).

\subsubsection*{(iii) Thin Film Superconductors in a Perpendicular External Field}
Here we consider an infinitely thin film superconductor occupying the
two-dimensional domain $\Omega$
in the $x_3=0$ plane. With $b_e(t)$
the normal to the film component of the given non-stationary
uniform external
magnetic field,
the normal to the film component of the total magnetic
field can then be expressed by the Biot--Savart law as
\begin{align}
b_3(\vx,t) = 
b_e(t) + \frac{1}{4\pi} \, \nabla \times \int_{\Omega}
\frac{\vj (\vy,t)}{|\vx-\vy|}
\, \dy\,,
\label{BS}
\end{align}
where $\vj$ is the sheet current density in the film.
Using Faraday's law with $\ve$, the component of the electric field tangential to the film,
and the change of variable $\vq$ in (\ref{curldiv}), we obtain that
\begin{align}
\frac{\partial b_3}{\partial t} = - \nabla \times \ve = - \del\,.\,\vq\,.
\label{Fltf}
\end{align}
As $\del\,.\,\vj =0$ in $\Omega$, which is simply connected,
we can introduce a stream (magnetization) function $w$,
which vanishes on $\partial \Omega$,
such that $\vj= \del \times w$ in $\Omega$.
Substituting this and (\ref{Fltf}) into the time derivative of (\ref{BS}), we obtain that
\begin{align}
 \frac{1}{4\pi}\, \nabla \times \int_{\Omega}
\frac{1}{|\vx-\vy|}
\, \del \times \frac{\partial w(\vy,t)}{\partial t}  \dy
+ \del\,.\,\vq(\vx,t) = - \frac{db_e(t)}{dt} \,.
\label{BStw}
\end{align}
The critical state model relations are given, as before, by (\ref{crit}).
However, in this problem we limit our considerations to the variational inequality case
and assume the Bean model with a field independent sheet critical current density
$j_c = k \in C(\overline{\Omega})$ and, as in (ii) above, $k(\vx) \geq k_{\rm min} >0$
for all $\vx \in \overline{\Omega}$.
The model relations can be reformulated as (\ref{cond3}) with $M(w)$ replaced by $k$,
and this supplemented with (\ref{BStw}) yields the mixed formulation
of this thin film superconductor problem.
For the initial data, we take $w(\cdot,0)=w_0(\cdot)$ with $|\del w_0| \leq k$.
Similarly, one can show that $w$ satisfies a primal variational inequality
problem, see Theorem \ref{conthmtfM} below.
In addition, one can approximate the critical state model relations
by the power law model (\ref{power}), see Barrett and Prigozhin
\cite{EFILMMATH} for further details and subsection \ref{Outline} below.
Similarly to \cite{BP5}, we note that the sign of $\vq$ is changed in
\cite{EFILMMATH} ($\vv$ in the notation there).

\subsection{Notation}

Above, and throughout, we adopt the standard notation for Sobolev spaces
on a bounded domain $D \subset {\mathbb R}^d$ with a Lipschitz boundary,
denoting the norm of
$W^{m,s}(D)$ ($m \in {\mathbb N}$, $s\in [1, \infty]$) by
$\|. \|_{m,s,D}$ and the semi-norm by $|\cdot |_{m,s,D}$.
Of course, we have that
$|\cdot |_{0,s,D} \equiv \|\cdot \|_{0,s,D}$.
We extend these norms and semi-norms in the natural way to the corresponding
spaces of vector
functions.
For $s=2$, $W^{m,2}(D)$ will be denoted by
$H^m(D)$ with the associated norm and semi-norm written
as, respectively, $\|\cdot\|_{m,D}$ and $|\cdot|_{m,D}$.
We set $W^{1,s}_0(D):= \{\eta \in W^{1,s}(D) : \eta = 0 \mbox{ on } \partial D \}$,
and $H^1_0(D) \equiv W^{1,2}_0(D)$.
We recall the Poincar\'{e} inequality for any $s \in [1,\infty]$
\begin{align}
|\eta|_{0,s,D} \leq C_\star(D)\,|\del \eta|_{0,s,D}\qquad \forall \ \eta \in W^{1,s}_0(D)\,,
\label{Poin}
\end{align}
where the constant $C_\star(D)$ depends on $D$, but is independent of $s$;
see e.g.\  page 164 in Gilbarg and Trudinger \cite{GT}.
In addition, $|D|$ will denote the measure of $D$.
We require also $H^{\frac{1}{2}}(D)$ for $D \subset {\mathbb R}^2$ and
\begin{align}
H^{\frac{1}{2}}_{00}(D) := \left\{ \eta \in H^{\frac{1}{2}}(D) :
\widetilde{\eta} := \left\{ \begin{array}{l}
\eta \mbox{ in } D \\
0 \mbox{ in } {\mathbb R}^2 \setminus D
\end{array}
\right.
 \in H^{\frac{1}{2}}({\mathbb R}^2)
\right\}.
\label{Hhalf}
\end{align}
For any Banach space ${\cal B}$, we denote its dual by ${\cal B}^\star$.
Then we recall that
\begin{align}
\|\eta\|_{[H^{\frac{1}{2}}(D)]^\star}
\leq \left[ \|\eta\|_{[H^{1}(D)]^\star}\,\|\eta\|_{L^2(D)} \right]^{\frac{1}{2}}
\qquad \forall \ \eta \in L^2(D)\,.
\label{Hhalfstarnorm}
\end{align}

For $m\in \mathbb N$, let
(i) $C^m(\overline{D})$
denote the Banach space of continuous functions with
all derivatives up to order $m$
continuous on $\overline{D}$,
(ii) $C^m_0(D)$ denote the space of continuous functions with
compact support in $D$ with all derivatives up to order $m$
continuous on $D$ and (iii)
$C^m_0(\overline{D})$ denote
the Banach space 
$ \{\eta \in C^m(\overline{D}) : \eta = 0 \mbox{ on } \partial D\}$.
In the case $m=0$, we drop the superscript $0$ for all three spaces.

As one can identify $L^{1}(D)$ as a closed subspace of
the Banach space of bounded Radon measures,
${\cal M}(\overline{D}) \equiv [C(\overline{D})]^\star$,
it is convenient to adopt the notation
\begin{align}
\int_{\overline{D}} |\mu| \equiv \|\mu\|_{{\cal M}(\overline{D})}
:= \sup_{\stackrel{\eta \in C(\overline{D})}{|\eta|_{0,\infty,D}\leq 1}}
\langle \mu,\eta \rangle_{C(\overline{D})}
< \infty,
\label{Mnorm2}
\end{align}
where $\langle \cdot,\cdot \rangle_{{\cal B}}$
denotes the duality pairing on ${\cal B}^\star
\times {\cal B}$ for any Banach space ${\cal B}$.
We note that if
$\{\mu_n\}_{n \geq 0}$ is a bounded sequence in ${\cal M}(\overline{D})$,
then there exist a subsequence $\{\mu_{n_j}\}_{n_j \geq 0}$ and a
$\mu \in {\cal M}(\overline{D})$ such that as $n_j \rightarrow \infty$
\begin{align}
\mu_{n_j} \rightarrow \mu \quad \mbox{weakly in } {\cal M}(\overline{D});
\quad \mbox{i.e.}
\quad \langle \mu_{n_j}-\mu,\eta\rangle_{C(\overline{D})} \rightarrow 0\quad
\forall \ \eta \in C(\overline{D})\,.
\label{Mweak1}
\end{align}
In addition,
we have that
\begin{align}
\liminf_{n_j \rightarrow \infty} \int_{\overline{D}} |\mu_{n_j}|
\geq \int_{\overline{D}} |\mu|\,;
\label{Mweak2}
\end{align}
see e.g.\ page 223 in
Folland \cite{Folland}.
For $D \subset {\mathbb R}^2$ we require the following Banach spaces
\eqlabstart
\begin{alignat}{2}
\vV^s(D)&:= \{ \vv \in [L^s(D)]^2 : \del \,.\,\vv \in L^2(D)
\} \qquad &&\mbox{for a given } s \in [1,\infty]
\label{vVs}\\
\mbox{and} \qquad \qquad
\vV^{\cal M}(D)&:= \{ \vv \in [{\cal M}(\overline{D})]^2 : \del \,.\,\vv \in L^2(D)\}\,;
\label{vVM}\\
\vZ^s(D)&:= \{ \vv \in [L^s(D)]^2 : \del \,.\,\vv \in [H^{\frac{1}{2}}_{00}(D)]^\star
\} \qquad &&\mbox{for a given } s \in [1,\infty]
\label{vZs}\\
\mbox{and} \qquad \qquad
\vZ^{\cal M}(D)&:= \{ \vv \in [{\cal M}(\overline{D})]^2 : \del \,.\,\vv \in \in [H^{\frac{1}{2}}_{00}(D)]^\star\}\,.
\label{vZM}
\end{alignat}
\eqlabend

We recall 
the Aubin--Lions--Simon compactness theorem, see Corollary 4 in Simon \cite{Simon}.
Let ${\cal B}_0$, ${\cal B}$ and
${\cal B}_1$ be Banach
spaces, ${\cal B}_i$, $i=0,1$, reflexive, with a compact embedding ${\cal B}_0
\hookrightarrow {\cal B}$ and a continuous embedding ${\cal B} \hookrightarrow
{\cal B}_1$. Then, for $\alpha>1$, the embedding
\begin{align}
&\{\,\eta \in L^{\infty}(0,T;{\cal B}_0): \frac{\partial \eta}{\partial t}
\in L^{\alpha}(0,T;{\cal B}_1)\,\} \hookrightarrow C([0,T];{\cal B})
\label{compact1}
\end{align}
is compact.
We write $(\cdot,\cdot)$ for
the standard inner product on $L^2(\Omega)$.
Finally, throughout $C$ denotes a generic positive constant independent of
the power parameters, $r \in (1,2)$ and $p \in (2,\infty)$, recall (\ref{power}),
the mesh parameter $h$ and the time step
parameter $\tau$. Whereas, $C(s)$ denotes a positive constant dependent on
the parameter $s$.

\subsection{Outline}
\label{Outline}

We introduce
\eqlabstart
\begin{align}
c(\vpsi,\vv)&:=
\frac{1}{4\pi}\int_{\Omega}\int_{\Omega}
\frac{\vpsi(\vx)\,.\, \vv(\vy)}{|\vx-\vy|}
\,\dx\,\dy \label{bilformc}\\
\mbox{and} \qquad
a(\phi,\eta) &:= c(\del \,\phi,\del\,\eta) =
\frac{1}{4\pi}\int_{\Omega}\int_{\Omega}
\frac{\del\,\phi(\vx)\,.\, \del\,\eta(\vy)}{|\vx-\vy|}
\,\dx\,\dy\,.
\label{bilforma}
\end{align}
\eqlabend
It follows that $c(\cdot\cdot)$ and $a(\cdot,\cdot)$
are symmetric, continuous and coercive bilinear
forms on $[\,[H^{\frac{1}{2}}(\Omega)]^\star]^2\times [\,[H^{\frac{1}{2}}(\Omega)]^\star]^2$ and
$H^{\frac{1}{2}}_{00}(\Omega)\times H^{\frac{1}{2}}_{00}(\Omega)$, respectively,
see \cite[Lemma 2.1]{NA}. 
Then we introduce for all $\chi,\,\eta \in W^{1,p}_0(\Omega)$
\begin{align}
{\cal A}(\chi,\eta) := \left\{\begin{array}{ll}
(\chi,\eta) &\qquad \mbox{in cases (i) and (ii),}\\[2mm]
a\,(\chi,\eta) &\qquad \mbox{in case (iii);}
\end{array}
\right.
\label{calA}
\end{align}
and set $\|\cdot\|_{{\cal A}} = [{\cal A}(\cdot,\cdot)]^{\frac{1}{2}}$.
In addition,
we introduce for all $\eta \in L^2(0,T;C(\overline{\Omega}))$
\begin{align}
{\mathfrak M} (\eta)(\vx,t)\: := \left\{\begin{array}{ll}
M_\eps(\eta(\cdot,t))(\vx) &\qquad \mbox{in case (i),}\\[2mm]
k(\vx)\,\widehat M(\eta(\vx,t)+b_e(t)) &\qquad \mbox{in case (ii),}\\[2mm]
k(\vx)&\qquad \mbox{in case (iii)}
\end{array}
\right.
\qquad \mbox{for a.e.\ } (\vx,t) \in \Omega \times (0,T);
\label{calMt}
\end{align}
where $M_\eps(\cdot)$ is given by (\ref{Meps}), $\hat M(\cdot)$ satisfies (\ref{jcsc})
and $k \in L^\infty(\Omega)$ with $k(\vx) \geq k_{\rm min} >0$ for a.e.\ $\vx \in \Omega$.
We note that the assumption $\hat M_0 \in {\mathbb R}_{>0}$
does allow for any continuous $\hat{M}(\cdot)$ that is strictly 
positive on any bounded interval of $\mathbb{R}$, but such that
$\hat{M}(s) \rightarrow 0$ as $|s| \rightarrow \infty$.    
This follows as any solution of the critical state model
will be bounded, and hence $\hat M(\cdot)$ can be modified
to satisfy $\hat M_0 \in {\mathbb R}_{>0}$ without changing the problem;
see \cite{BP5} for details. 

Furthermore, we set
\eqlabstart
\begin{align}
w^0 &:= \left\{\begin{array}{ll}
w_0^\epsilon \in C^1_0(\overline{\Omega}) &\quad \mbox{in case (i),}\\[2mm]
w_0 \in W^{1,\infty}_0(\Omega) \quad \mbox{s.t.} \quad |\del w_0(\cdot)| \leq
{\mathfrak M}(w_0)(\cdot,0)&\quad \mbox{in cases (ii) and (iii);}
\end{array}
\right. \label{inidata}\\[2mm]
\qquad \mbox{and} \qquad
{\cal F} &:= \left\{\begin{array}{ll}
\mbox{nonnegative } f \in L^2(0,T;L^2(\Omega))&\quad \mbox{in case (i),}\\[2mm]
- \frac{d{b_e}}{dt} \in L^2(0,T) &\quad \mbox{in cases (ii) and (iii).}
\end{array}
\right.
\label{calFn}
\end{align}
\eqlabend
It follows from (\ref{bal}), (\ref{BStw}), (\ref{power}), (\ref{calA}), (\ref{bilformc},b),
(\ref{calMt}) and (\ref{inidata},b) that the formal weak mixed formulation of the
power law approximation of our three quasi-variational inequality problems
can be written in a unified way for a given $r \in (1,2)$:

(Q$_r$) Find $w_r \in H^{1}(0,T;W^{1,p}_0(\Omega))$ and $\vq_r \in L^r(0,T;[L^r(\Omega)]^2)$
such that for a.a.\ $t \in (0,T)$
\eqlabstart
\begin{alignat}{2}
{\cal A}(\frac{\partial w_r}{\partial t},\eta) - (\vq_r, \del \eta) & =({\cal F}, \eta)
\qquad &&\forall \ \eta \in W^{1,p}_0(\Omega),
\label{Qrv1}
\\
({\mathfrak M}(w_r)\,|\vq_r|^{r-2} \vq_r , \vv) + (\del w_r,\vv) &= 0
\qquad && \forall \ \vv \in [L^r(\Omega)]^2;
\label{Qrv2}
\end{alignat}
\eqlabend
where $w_r(\cdot,0) = w^0(\cdot)$.

We will be more precise about the function spaces of this weak formulation
with respect to the different problems (i), (ii) and (iii) in Section \ref{conv}.
In \cite{sandqvi}, \cite{BP5} and \cite{EFILMMATH}, we introduced
a finite element approximation of (\ref{Qrv1},b) based on the lowest order Raviart--Thomas
element for $\vq_r$ for cases (i), (ii) and (iii), respectively;
with piecewise constants for $w_r$ in cases (i) and (ii), and continuous piecewise
linears in case (iii).
There integration by parts was performed on the second terms on the left-hand sides of
(\ref{Qrv1},b) as the Raviart--Thomas element is a conforming approximation of the
divergence operator. In addition, in case (ii) we exploited (\ref{jcsc}) and
based our finite element approximation on the following rewrite of (\ref{Qrv2})
\begin{align}
(k\,|\vq_r|^{r-2} \vq_r , \vv) - (\hat F(w_r+b_e)-\hat F(b_e),\del\,.\,\vv) &= 0
\qquad && \forall \ \vv \in \vV^r(\Omega).
\label{Qrv2a}
\end{align}
In \cite{sandqvi} and \cite{BP5} we proved (subsequence) convergence
of these finite element approximations in cases (i) and (ii), respectively,
to the corresponding weak
mixed formulation of the critical state model, (Q), as the mesh parameters
tend to zero and $r \rightarrow 1$.
In \cite{EFILMMATH} we
proved convergence
of the finite element approximation in case (iii)
to the corresponding weak
mixed formulation of the power law model, (Q$_r$), as the mesh parameters
tend to zero. 
We note that in case (iii), one can show that the solution of (Q$_r$)
is unique as ${\mathfrak M}$ only depends on $\vx$, recall (\ref{calMt}).
We also proved in \cite{EFILMMATH} (subsequence) convergence of the solution to (Q$_r$) 
to a solution of the corresponding
weak mixed formulation of the critical state model, (Q), as $r \rightarrow 1$.
Finally, we remark that the power law model, (Q$_r$), is of interest in its own right in
the superconductivity context, cases (ii) and (iii), as it is a popular choice among engineers
for a current-voltage relation for some superconducting materials.

In this paper we consider a simpler finite element approximation of (\ref{Qrv1},b)
based on a nonconforming linear approximation of $w_r$ and a piecewise constant
approximation of $\vq_r$. Of course, for linear second order elliptic problems
the nonconforming linear approximation is a computationally inexpensive way of obtaining
the lowest order Raviart--Thomas approximation, see Marini \cite{Marini};
but this does not carry across to nonlinear problems.
We note that in \cite{sandqvi} for case (i), in addition to considering the Raviart--Thomas
approximation of (Q$_r$), (\ref{Qrv1},b), we also considered an approximation
based on continuous piecewise linears for $w_r$ and a piecewise constant
approximation of $\vq_r$. Once again, we showed (subsequence) convergence of this
finite element approximation to the corresponding weak mixed formulation of the sandpile model, (Q),
as the mesh parameters tend to zero and $r \rightarrow 1$.
Although this finite element approximation leads to a good approximation of the surface
$w$ in practice, the approximation of the sand flux $\vq$ is poor. We note that
all the convergence results stated above for the sand flux (rotated electric field) variable
are weak convergence results.
Hence there is no guarantee that this flux approximation will be useful in practice.
Nevertheless, the Raviart--Thomas sand flux (rotated electric field) approximations for
(i), (ii) and (iii) converged strongly in practice for the numerical experiments
in \cite{sandqvi}, \cite{BP5} and \cite{EFILMMATH}, respectively;
see also \cite{SUST} for case (iii).
Similarly, strong convergence is also observed
in practice for the sand flux (rotated electric field) approximation
resulting from the nonconforming linear approximation of $w_r$
and constant approximation $\vq_r$ studied in this paper.
For case (iii), see also \cite{Transport} where thin film problems
involving transport currents, which lead to  non-homogenous time-dependent boundary data for $w_r$
and singular time-dependent forcing data ${\cal F}$, are solved using this nonconforming
approximation.

The outline of this paper is as follows.
In the next section we introduce our nonconforming linear finite element
approximation, (Q$^{h,\tau}_r$), of the power law mixed formulation (Q$_r$), (\ref{Qrv1},b),
and prove well-posedness and stability bounds.
Here $h$ and $\tau$ are the spatial and temporal discretization
parameters, respectively.
In Section \ref{conv} we first prove (subsequence) convergence of (Q$^{h,\tau}_r$)
to (Q$^{\tau}_r$), a discrete time approximation of (Q$_r$),  as $h \rightarrow 0$.
Then under various assumptions, and appealing to results in
\cite{sandqvi}, \cite{BP5} and \cite{EFILMMATH} as much as possible,
we prove (subsequence) convergence of (Q$^\tau_r$) to (Q),
as $\tau \rightarrow 0$ and $r \rightarrow 1$,
for case (i); and (subsequence) convergence of (Q$^\tau_r$) to (Q$_r$),
as $\tau \rightarrow 0$,
and then (subsequence) convergence of (Q$_r$) to (Q),
as $r \rightarrow 1$, in cases (ii) and (iii).
The full sequence converges in case (iii) in the first two convergence  results, as in this case
one can prove uniqueness of the solution to problems (Q$^\tau_r$) and (Q$_r$).
Finally, in Section \ref{numexpts} we state an algorithm for solving the resulting
nonlinear algebraic equations arising from the approximation (Q$^{h,\tau}_r$) at each time level,
and present some numerical experiments.

\section{Finite Element Approximation}
\label{fea}
\setcounter{equation}{0}

We make the following assumptions on the data.
\vspace{2mm}

\noindent
{\bf (A1)} $\Omega \subset {\mathbb R}^2$
is simply connected and has a
Lipschitz boundary $\partial \Omega$
with outward unit normal $\underline{\nu}$.
The conditions stated on the data in (\ref{inidata},b) and (\ref{calMt}) hold.
In addition, in case (i)
the initial data $w_0^\epsilon \in 
C^1_0(\overline{\Omega})$ is such that
$\del w_0^\epsilon \,.\, \underline{\nu} < k_0$.

\vspace{2mm}

For ease of exposition, we shall assume that $\Omega$ is a polygonal domain to
avoid perturbation of domain errors in the finite element approximation. We make
the following standard assumption on the partitioning.
\vspace{2mm}

\noindent
{\bf (A2)} $\Omega$ is polygonal.
Let
 $\{{\cal T}^h\}_{h>0}$ be a regular family
 of partitionings of $\Omega$
 into disjoint open
 triangles $\sigma$ with $h_{\sigma}:={\rm diam}(\sigma)$
 and $h:=\max_{\sigma \in {\cal T}^h}h_{\sigma}$, so that
 $\overline{\Omega}=\cup_{\sigma\in{\cal T}^h}\overline{\sigma}$.
Moreover, $k \mid_{\sigma}$ can be extended to $k \in C(\overline{\sigma})$ for all
$\sigma \in {\cal T}^h$; that is, $k$ is piecewise continuous and its discontinuities only occur
along the internal edges of ${\cal T}^h$.

\vspace{2mm}
Let $\underline{\nu}_{\partial \sigma}$ be the outward unit normal to $\partial \sigma$,
the boundary of $\sigma$.
We then introduce the following finite element spaces
\eqlabstart
\begin{align}
S^h
&:= \{ \eta^h \in L^\infty(\Omega) : \eta^h \mid_{\sigma} = a_{\sigma} \in
{\mathbb R} \quad \forall \ \sigma \in {\cal T}^h \,\}\,,
\qquad \vS^h := [S^h]^2\,,
\label{Shc}\\
{U}^h
&:= \{ \eta^h \in C(\overline{\Omega}) : \eta^h \mid_{\sigma}
= \va_{\sigma}\,.\,\vx + b_{\sigma},
\ \va_{\sigma} \in {\mathbb R}^2, \ b_{\sigma} \in {\mathbb R}
\quad \forall \ \sigma \in {\cal T}^h \}\,,
\qquad U^h_0 := U^h \cap W^{1,\infty}_0(\Omega)\,,
\label{Wh}\\
{N}^h
&:= \{ \eta^h \in L^\infty(\Omega) : \eta^h \mid_{\sigma}
= \va_{\sigma}\,.\,\vx + b_{\sigma},
\ \va_{\sigma} \in {\mathbb R}^2, \ b_{\sigma} \in {\mathbb R}
\quad \forall \ \sigma \in {\cal T}^h, \nonumber \\
& \hspace{0.48in}
\mbox{ and } \eta^h
\mbox{ is continuous at the midpoints of the edges of neighbouring triangles}\}\,,
\label{Nh}\\
N^h_0 &:= \{ \eta^h \in N^h : \eta^h =0 \mbox{ at the midpoints of
the edges on } \partial \Omega\}\,. \label{Nh0}
\end{align}
\eqlabend

Let $\pi^h_U : C(\overline{\Omega})  \rightarrow U^h$ denote the $U^h$ interpolation
operator such that $\pi^h_U \eta(\vx_j^v)= \eta(\vx_j^v)$, $j=1,\ldots, J^v$,
where $\{\vx_j^v\}_{j=1}^{J^v}$ are the vertices of the partitioning ${\cal T}^h$.
Let $\pi^h_N : C(\overline{\Omega})  \rightarrow N^h$ denote the $N^h$ interpolation
operator such that $\pi^h_N \eta(\vx_j^e)= \eta(\vx_j^e)$, $j=1,\ldots, J^e$,
where $\{\vx_j^e\}_{j=1}^{J^e}$ are the midpoints of the edges of the partitioning ${\cal T}^h$.
We note for $m=0$ and 1 and any $s \in [1,2]$ that
\eqlabstart
\begin{alignat}{2}
&|(I-\pi^h_U) \eta|_{m,\sigma}
+|(I-\pi^h_N) \eta|_{m,\sigma}
\leq C\,h_\sigma^{3-m-\frac{2}{s}}\,|\eta|_{2,s,\sigma}
\qquad &&\forall \ \sigma \in {\cal T}^h, 
\label{piconv0}\\
&\lim_{h \rightarrow 0} \left[
\|(I-\pi^h_U) \eta\|_{m,\infty,\Omega}
+
|(I-\pi^h_N) \eta|_{0,\infty,\Omega}
+ m\,|\del \eta - \del_h (\pi^h_N \eta)|_{0,\infty,\Omega}
\right]
= 0
\qquad &&\forall \ \eta \in C^m(\overline{\Omega})
\,;
\label{piconv}
\end{alignat}
\eqlabend
where $I$ is the identity operator and
\begin{align}
\del_h \eta^h \mid_{\sigma} = \del \eta^h \qquad \forall \ \sigma \in {\cal T}^h,
\qquad \forall \ \eta^h \in N^h\,.
\label{delh}
\end{align}

Let $\vP^h: [L^1(\Omega)]^2 \rightarrow \vS^h$ be such that
\begin{align}
\vP^h \vv \mid_{\sigma} = \mints \vv \qquad
\forall \ \sigma \in {\cal T}^h\,,
\label{Ph}
\end{align}
where $\mintD \,\cdot :=
\frac{1}{|D|} \int_{D} \cdot \,{\rm d}\vx.$
We note that
\eqlabstart
\begin{align}
&|\vP^h \vv|_{0,s,\sigma} \leq |\vv|_{0,s,\sigma} \qquad \forall\ \vv \in [L^s(\sigma)]^2,
\quad s \in [1,\infty],\quad \forall \ \sigma \in {\cal T}^h\,,
\label{Phstab}\\
&\lim_{h \rightarrow 0}
|\,|\vv|-|\vP^h\vv|\,|_{0,\infty,\Omega}
\leq
\lim_{h \rightarrow 0}
|\vv-\vP^h \vv|_{0,\infty,\Omega} = 0
\qquad \forall \ \vv \in [C(\overline{\Omega})]^2\,.
\label{Pconv}
\end{align}
\eqlabend
In addition, one can show by mapping to a reference element, applying a trace inequality
and the Poincar\'{e} inequality (\ref{Poin}), and then mapping back that for any $s \in [1,\infty]$
and for all $\sigma \in {\cal T}^h$
\begin{align}
|(\vI-\vP^h)\vv|_{0,s,\partial_i \sigma} \leq C \,h_\sigma^{1-\frac{1}{s}}\,
|\vv|_{1,s,\sigma}
\qquad \forall \ \vv \in [W^{1,s}(\sigma)]^2,\qquad i=1,\,2,\,3\,,
\label{Phpis}
\end{align}
where $\partial_i \sigma$ is one of the three edges of $\partial \sigma$; that is
$\partial \sigma = \sum_{i=1}^3 \partial_i \sigma$.
Similarly, we define $P^h : L^1(\Omega) \rightarrow S^h$ with the equivalent
to (\ref{Phstab},b) and (\ref{Phpis}) holding.
In addition, we have that for any $s \in [1,\infty]$ and for all $\sigma \in {\cal T}^h$
\begin{align}
|\,(I-\mintpis)\,\eta^h|_{0,s,\partial_i \sigma} &\leq C\, h_\sigma\,
|\del \eta^h|_{0,s,\partial \sigma}
\leq C\, h_\sigma^{1-\frac{1}{s}}\,
|\del \eta^h|_{0,s,\sigma}  \qquad \forall \ \eta^h \in N^h ,\qquad i=1,\,2,\,3\,.
\label{mintsetahbd}
\end{align}

We recall for $r>1$ and
for all $\vc, \,\vd \in {\mathbb R}^d$ that
\eqlabstart
\begin{alignat}{2}
\frac{1}{r}\,\frac{\partial|\vc|^r}{\partial c_i}  = |\vc|^{r-2}\,c_i
\qquad &\Rightarrow
\qquad
&&|\vc|^{r-2}\,\vc\,
.\,(\vc-\vd) \geq \ts \frac{1}{r}\,\left[\, |\vc|^r -|\vd|^r\,\right]
\geq
|\vd|^{r-2}\,\vd\,
.\,(\vc-\vd)\,,
\label{aederiv}
\\
\qquad &\mbox{and}\qquad
&&\left( |\vc|^{r-2}\,\vc - |\vd|^{r-2}\,\vd \right)\,.\,\vc
\geq
\ts \frac{r-1}{r}
\,
[\,|\vc|^r - |\vd|^r\,]\,.
\label{abcoer}
\end{alignat}
\eqlabend
Let
$0= t_0 < t_1 < \ldots < t_{N-1} < t_N = T$ be a
partitioning of $[0,T]$
into possibly variable time steps $\tau_n := t_n -
t_{n-1}$, $n=1,\ldots, N$. We set
$\tau := \max_{n=1,\ldots, N}\tau_n$ and. on recalling (\ref{calFn}), we introduce
\begin{align}
{\cal F}^n(\cdot) := \frac{1}{\tau_n} \int_{t_{n-1}}^{t_n} {\cal F}(\cdot,t) \,{\rm d}t
\in L^2(\Omega) \qquad n=1,\ldots, N\,.
\label{fn}
\end{align}
We note  that
\begin{align}
\sum_{n=1}^N \tau_n \,|{\cal F}^n|_{0,s,\Omega}^s
\leq \int_0^T |{\cal F}|_{0,s,\Omega}^s\,{\rm d}t
\qquad \mbox{for any } s \in [1,2]\,.
\label{fnsum}
\end{align}

On setting
\begin{align}
w^{\eps,h}_0 = P^h [\pi^h_N w^\eps_0]\,,
\label{weps0h}
\end{align}
we introduce $M^h_\eps : S^h \rightarrow S^h$ approximating $M_\eps : C(\overline{\Omega})
\rightarrow C(\overline{\Omega})$, defined by (\ref{Meps}), for any $\sigma \in {\cal T}^h$ as
\begin{align}
&M_{\varepsilon}^h(\eta^h) \mid_{\sigma}
:=\left\{ \begin{array}{ll}  k_0 \quad& \eta^h_\sigma \geq w_{0,\sigma}^{\epsilon,h}
+ \varepsilon,\\
k_{1,\sigma}^{\epsilon,h} + (k_0-k_{1,\sigma}^{\epsilon,h})
\left(\ds\frac{\eta^h_\sigma-w_{0,\sigma}^{\epsilon,h}}{\varepsilon}\right)
\quad & \eta^h_\sigma \in [w_{0,\sigma}^{\epsilon,h},w_{0,\sigma}^{\epsilon,h}+\varepsilon],\\[4mm]
k_{1,\sigma}^{\epsilon,h}
:= \max(k_0,|\del_h (\pi^h_N w_0^\epsilon) \mid_{\sigma}\!|)
\quad & \eta^h_\sigma\le w_{0,\sigma}^{\epsilon,h};\end{array}\right.
\label{Meps_h}
\end{align}
where $\eta^h_\sigma = \eta^h \mid_{\sigma}$ and $w^{\epsilon,h}_{0,\sigma} =
 w^{\epsilon,h}_{0} \mid_{\sigma}$ for all $\sigma \in {\cal T}^h$.

We note that $M_\eps$ is also well-defined on $S^h$ with $M_\eps : S^h \rightarrow L^\infty(\Omega)$,
and we have the following result.
\begin{lmm}\label{Mepshapproxlem}
For any $\eta^h \in S^h$, we have  that
\begin{align}
|M_\eps(\eta^h) - M^h_\eps(\eta^h)|_{0,\infty,\Omega} 
\leq C(\epsilon^{-1}) \left[\, |(I-P^h) w^{\eps}_0|_{0,\infty,\Omega}+
|\del w^{\eps}_0 - \del_h(\pi^h_N w^{\eps}_0)|_{0,\infty,\Omega}\,
\right].
\label{Mepshapprox}
\end{align}
\end{lmm}
\begin{proof}
See the proof of Lemma 2.1 in Barrett and Prigozhin \cite{sandqvi}.
\end{proof}

Finally, it follows from (\ref{Meps_h}), (\ref{piconv}) and Assumption (A1) that
\begin{align}
k^{\epsilon,h}_{1,\infty} := \max_{\sigma \in {\cal T}^h} k^{\epsilon,h}_{1,\sigma} \leq C\,.
\label{kepsh1max}
\end{align}

\subsection{Approximation (Q$^{h,\tau}_{r}$)}

On recalling (\ref{calA}) and (\ref{bilformc},b), we introduce for all $\chi^h,\,\eta^h \in N^h_0$
\begin{align}
{\cal A}^h (\chi^h,\eta^h) := \left\{\begin{array}{ll}
(\chi^h,\eta^h) &\qquad \mbox{in cases (i) and (ii),}\\[2mm]
c\,(\del_h \chi^h,\del_h \eta^h) &\qquad \mbox{in case (iii);}
\end{array}
\right.
\label{calAh}
\end{align}
and set $\|\cdot\|_{{\cal A}^h} = [{\cal A}^h(\cdot,\cdot)]^{\frac{1}{2}}$.
In addition, on recalling (\ref{calMt}), we introduce for all $\eta^h \in S^h$
\begin{align}
{\mathfrak M}^{h,n} (\eta^h) \mid_{\sigma}\: := \left\{\begin{array}{ll}
M^h_\eps(\eta^h) \mid_\sigma &\qquad \mbox{in case (i),}\\[2mm]
k(\vx_\sigma)\,\widehat M(\eta^h \mid_\sigma +b_e(t_n)) &\qquad \mbox{in case (ii),}\\[2mm]
k(\vx_\sigma)&\qquad \mbox{in case (iii)}
\end{array}
\right.
\qquad \forall \ \sigma \in {\cal T}^h,
\label{calMh}
\end{align}
where $\vx_\sigma$ is the centroid of $\sigma$.
We note from (\ref{calMh}), (\ref{Meps_h}), (\ref{kepsh1max}) and Assumption (A1) that
${\mathfrak M}^{h,n} : S^h \rightarrow S^h$ and
there exist ${\mathfrak M}_{\rm min},\,{\mathfrak M}_{\rm max} \in {\mathbb R}$ such that
for $n=1,\ldots,N$
\begin{align}
0 < {\mathfrak M}_{\rm min} \leq {\mathfrak M}^{h,n}(\eta^h) \mid_{\sigma} \,
\leq {\mathfrak M}_{\rm max}
\qquad \forall \ \eta^h \in S^h, \qquad \forall \ \sigma \in {\cal T}^h\,.
\label{calMhb}
\end{align}

We now define our finite element approximation of (Q$_r$), (\ref{Qrv1},b), for a given $r>1$:

\noindent
{\bf (Q$^{h,\tau}_{r}$)}
For $n = 1,\ldots, N$, find $W^n_{r} \in N^h_0$ and
$\vQ^n_{r} \in \vS^h$ such that
\eqlabstart
\begin{alignat}{2}
{\cal A}^h\left(\frac{W^n_{r}-W^{n-1}_{r}}{\tau_n},\eta^h\right)
- (\vQ^n_{r},\del_h\eta^h) &= ({\cal F}^n
,\eta^h)
\quad \;&&\forall \ \eta^h \in N^h_0,\label{Qthr1B}\\
({\mathfrak M}^{h,n}(P^h W^n_{r})\,|\vQ^n_{r}|^{r-2} \vQ^n_{r},
\vv^h) + (\del_h W^n_{r},\vv^h)  &=0
\quad \;
&&\forall \ \vv^h\in \vS^h\,;
\label{Qthr2B}
\end{alignat}
\eqlabend
where $W^0_r = \pi^h_N w^0$.

Associated with (Q$^{h,\tau}_{r}$) is the corresponding
approximation of a generalised $p$-Laplacian
problem for $p>1$,
where we recall that $\frac{1}{r}+ \frac{1}{p}=1$:

\noindent
{\bf (P$^{h,\tau}_{p}$)}
For $n = 1,\ldots, N$, find $W^n_{r} \in N^h_0$ such that
\begin{align}
&{\cal A}^h\left(\frac{W^n_{r}-W^{n-1}_{r}}{\tau_n},\eta^h\right)
+ \left( [{\mathfrak M}^{h,n}(P^h W^n_{r})]^{-(p-1)}\, |\del_h W^n_{r}|^{p-2} \,
\del_h W^n_{r}, \del_h \eta^h \right)
= ({\cal F}^n ,\eta^h)
\qquad \forall \ \eta^h \in N^h_0\,,\label{PthrB}
\end{align}
where $W^0_r = \pi^h_N w^0$.

\begin{thrm}\label{Qthrstab}
Let the Assumptions (A1) and (A2) hold.
Then for all $r \in (1,2)$,
for all regular partitionings ${\cal T}^h$ of $\Omega$,
and for all $\tau_n >0$,
there exists a solution,
$W^n_{r} \in N^h_0$ and $\vQ^n_{r} \in \vS^h$
to the $n^{\rm th}$ step
of (Q$^{h,\tau}_{r}$), (\ref{Qthr1B},b).
This solution is unique in case (iii).
In addition, we have that
\begin{align}
\max_{n=0,\ldots, N}
\|W^n_{r}\|_{{\cal A}^h} + \sum_{n=1}^N \|W^n_{r}-W^{n-1}_{r}\|_{{\cal A}^h}^2 +
\sum_{n=1}^N \tau_n\,|\vQ^n_{r}|_{0,r,\Omega}^r
+\left( \sum_{n=1}^N \tau_n |\del_h W^n_{r}|_{0,p,\Omega}^p\right)^{\frac{1}{p}}
&\leq C 
\,
\label{energy3}
\end{align}
where $\frac{1}{r}+ \frac{1}{p}=1$.
Moreover, 
(Q$^{h,\tau}_{r}$), (\ref{Qthr1B},b), is equivalent to (P$^{h,\tau}_{p}$), (\ref{PthrB}).
\end{thrm}
\begin{proof}
The proof is similar to the proof of Theorem 2.2 in \cite{sandqvi} with $U^h_0$
replaced by $N^h_0$.
It follows immediately from (\ref{Qthr2B}) that
for all $\sigma \in {\cal T}^h$
\begin{align}
\del W^n_{r} &=- {\mathfrak M}^{h,n}(P^h W^n_{r})\,|\vQ^n_{r}|^{r-2} \vQ^n_{r} 
\quad  \Leftrightarrow \quad \vQ^n_{r}
= -
[{\mathfrak M}^{h,n}(P^h W^n_{r})]^{-(p-1)}\, |\del W^n_{r}|^{p-2} \, \del W^n_{r}
\quad \mbox{ on } \sigma\,.
\label{QthrW}
\end{align}
Substituting this expression for $\vQ^n_{r}$ into (\ref{Qthr1B}) yields (\ref{PthrB}).
Hence (P$^{h,\tau}_{p}$), with (\ref{QthrW}), is equivalent to (Q$^{h,\tau}_{r}$).

Consider the strictly convex minimization
problem:
\eqlabstart
\begin{align}
\min_{\eta^h \in N^h_0} E^{h,n}_p(\eta^h)\,,
\label{pmin}
\end{align}
where, for a given $\varphi^h \in N^h_0$, $E^{h,n}_p : N^h_0 \rightarrow {\mathbb R}$ is defined by
\begin{align} E^{h,n}_p(\eta^h):=
\frac{1}{2 \tau_n}\, \|\eta^h-W^{n-1}_r\|_{{\cal A}^h}^2 + \frac{1}{p} \int_{\Omega}
[{\mathfrak M}^{h,n}(P^h \varphi^h)]^{-(p-1)}\,
|\del_h \eta^h|^p \,{\rm d}\vx -({\cal F}^n,\eta^h)\,.
\label{pminE}
\end{align}
\eqlabend
In case (iii), as, on recalling (\ref{calMh}),
${\mathfrak M}^{h,n}(\cdot)$ only depends on $\vx$,
(\ref{PthrB}) is the Euler-Lagrange system associated
with the strictly convex minimization problem (\ref{pmin},b).
Hence, in case (iii)
there exists a unique solution to  (P$^{h,\tau}_{p}$), (\ref{PthrB}),
and therefore to (Q$^{h,\tau}_{r}$), (\ref{Qthr1B},b).

We now apply the Brouwer fixed point theorem to prove existence of a solution to  (P$^{h,\tau}_{p}$),
and therefore to (Q$^{h,\tau}_{r}$) in cases (i) and (ii).
Let $F^h : N^h_0 \rightarrow N^h_0$ be such that for any $\varphi^h \in
N^h_0$, $F^h \varphi^h \in N^h_0$ solves
\begin{align}
&{\cal A}^h\left(\frac{F^h \varphi^h -W^{n-1}_{r}}{\tau_n},\eta^h\right)
+ \left( [{\mathfrak M}^{h,n}(P^h\varphi^h)]^{-(p-1)}\, |\del_h F^h \varphi^h |^{p-2} \,
\del_h  F^h \varphi^h, \del_h \eta^h \right)
= ({\cal F}^n ,\eta^h)
\qquad \forall \ \eta^h \in N^h_0\,.\label{Fh}
\end{align}
The well-posedness of the mapping $F^h$ follows from noting that
(\ref{Fh}) is the Euler--Lagrange system associated with
the strictly convex minimization problem (\ref{pmin},b).
that is, there exists a unique element $F^h \varphi^h \in N^h_0$ solving (\ref{Fh}).
It follows immediately from (\ref{pmin},b), as $\|\cdot\|_{{\cal A}^h}
\equiv |\cdot|_{0,\Omega}$ in cases (i) and (ii), that
\begin{align}
\frac{1}{2 \tau_n}\, |F^h \varphi^h-W^{n-1}_{r}|_{0,\Omega}^2  -({\cal F}^n,F^h \varphi^h)
\leq
E^{h,n}_p(F^h \varphi^h) \leq E^{h,n}_p(0) = \frac{1}{2 \tau_n}\, |W^{n-1}_{r}|_{0,\Omega}^2\,.
\label{pminEb}
\end{align}
It is easily deduced from (\ref{pminEb}) 
that
\begin{align}
F^h \varphi^h \in B_\gamma := \{\eta^h \in N^h_0 : |\eta^h|_{0,\Omega} \leq \gamma \}\,,
\label{Bgamma}
\end{align}
where $\gamma \in {\mathbb R}_{>0}$
depends on $|W^{n-1}_{r}|_{0,\Omega}$, $|{\cal F}^n|_{0,\Omega}$ and
$\tau_n$.
Hence $F^h : B_\gamma \rightarrow B_\gamma$. In addition, it is easily verified
that the mapping $F^h$ is continuous, as ${\mathfrak M}^{h,n}: S^h \rightarrow S^h$
is continuous with respect to $S^h$ on recalling (\ref{calMh}), (\ref{Meps_h}) and (\ref{jcsc}).
Therefore, the Brouwer fixed point theorem yields that the mapping $F^h$ has at least one
fixed point in $B_\gamma$. Hence,
there exists a solution to  (P$^{h,\tau}_{p}$), (\ref{PthrB}),
and therefore to (Q$^{h,\tau}_{r}$), (\ref{Qthr1B},b), in cases (i) and (ii).

It follows from (\ref{QthrW}) and (\ref{calMhb}) that for $n=1,\ldots, N$
\begin{align}
|\del_h W^n_{r}|_{0,p,\Omega}^p &=
|[{\mathfrak M}^{h,n}(P^h W^n_{r})]^{p-1}\vQ^n_{r}|_{0,r,\Omega}^r
\leq ({\mathfrak M}_{\rm max})^{p-1}\,
({\mathfrak M}^{h,n}(P^h W^n_{r}),|\vQ^n_{r}|^r)
\,.
\label{WnQn}
\end{align}
Choosing $\eta^h = W^n_{r}$, $\vv^h= \vQ^{n}_{r}$ in (\ref{Qthr1B},b),
combining and noting 
the simple
identity
\begin{align}
(c-d)\,c = \frac{1}{2} \left[ c^2 + (c-d)^2 - d^2 \right]
\qquad \forall \ c,\,d \in {\mathbb R}
\,,
\label{simpid}
\end{align}
we obtain for $n=1,\ldots, N$, on applying a Young's inequality and (\ref{Poin}), that
for all $\delta >0$
\begin{align}
&\|W^n_{r}\|_{{\cal A}^h}^2 + \|W^n_{r}-W^{n-1}_{r}\|_{{\cal A}^h}^2
+ 2\tau_n\, ({\mathfrak M}^{h,n}(P^h W^n_{r}),|\vQ^n_{r}|^r)
\nonumber \\
& \hspace{2in} = \|W^{n-1}_{r}\|_{{\cal A}^h}^2 + 2 \tau_n\,({\cal F}^n,W^n_{r})
\nonumber \\
& \hspace{2in} \leq \|W^{n-1}_{r}\|_{{\cal A}^h}^2 +
2 \tau_n\,\left[\,\frac{1}{r}\,\delta^{-r}\,|{\cal F}^n|_{0,r,\Omega}^r
+ \frac{1}{p}\,\delta^p |W^n_{r}|_{0,p,\Omega}^p\,\right]
\nonumber \\
& \hspace{2in} \leq \|W^{n-1}_{r}\|_{{\cal A}^h}^2
+ 2 \tau_n\,\left[\,\frac{1}{r}\,\delta^{-r}\,|{\cal F}^n|_{0,r,\Omega}^r
+ \frac{1}{p}\,[\delta\,C_\star(\Omega)]^p |\del W^n_{r}|_{0,p,\Omega}^p\,\right]\,.
\label{Wnr1}
\end{align}
It follows on summing (\ref{Wnr1}) from $n=1$ to $m$, with $\delta = 1/
(C_\star(\Omega)\,[{\mathfrak M}_{\rm max}]^{\frac{1}{r}})$,
and noting (\ref{WnQn}) and (\ref{calMhb}) that for $m=1,\ldots, N$
\begin{align}
&\|W^m_{r}\|_{{\cal A}^h}^2 + \sum_{n=1}^m \|W^n_{r}-W^{n-1}_{r}\|_{{\cal A}^h}^2
+ \sum_{n=1}^m
\tau_n\, ({\mathfrak M}^{h,n}(P^h W^n_{r}),|\vQ^n_{r}|^r)
\nonumber \\
& \hspace{3.2in}
\leq \|W^0_{r}\|_{{\cal A}^h}^2
+ 2\,[C_\star(\Omega)]^r\, {\mathfrak M}_{\rm max}
\sum_{n=1}^m \tau_n\,|{\cal F}^n|_{0,r,\Omega}^r\,.
\label{Wnr2}
\end{align}
The desired result (\ref{energy3}) follows immediately from (\ref{Wnr2}),
(\ref{fnsum}), (\ref{calMhb}) and   (\ref{WnQn}).
\end{proof}

We end this section with the following discrete
Poincar\'{e} and compactness results for $N^h_0$,
which are extensions of
Proposition 4.13 in Chapter 1 and
Theorem 2.4 in Chapter 2 of Temam \cite{TemamNS}.
In addition, we are more precise about the domain $\Omega$ and the subsequent
elliptic regularity.

\begin{lmm}\label{NhPoinlem}
Let $s \in (1,\infty)$ and the Assumption (A2) hold. Then we have that
\begin{align}
|\,(\eta^h, \del \,.\, \vv) + (\del_h \eta^h,\vv)|
\leq C\,h\,|\del_h \eta^h|_{0,s,\Omega}\,|\vv|_{1,s',\Omega}
\qquad \forall \ \eta^h \in N^h_0, \quad
\forall \
\vv \in [W^{1,s'}(\Omega)]^2\,;
\label{intbyparts}
\end{align}
where, here and throughout the paper, $\frac{1}{s}+\frac{1}{s'}=1$.
Hence, it follows that
\begin{align}
|\eta^h|_{0,s,\Omega} \leq C\,|\del_h \eta^h|_{0,s,\Omega} \qquad
\forall \ \eta^h \in N^h_0\,.
\label{NhPoin}
\end{align}
\end{lmm}
\begin{proof}
First on splitting $\partial \sigma$ into its three edges, i.e.\
$\partial \sigma = \sum_{i=1}^3 \partial_i \sigma$, it follows
from (\ref{mintsetahbd}) and (\ref{Phpis}) that
for all $\eta^h \in N^h_0$ and for all $\vv \in [W^{1,s'}_0(\Omega)]^2$
\begin{align}
\left|\,(\eta^h, \del \cdot \vv) + (\del_h \eta^h,\vv)\right|
&= \left|\sum_{\sigma \in {\cal T}^h} \sum_{i=1}^3 \int_{\partial_i \sigma} \eta^h\,
\,\vv \,.\,\underline{\nu}_{\partial_i \sigma} \,{\rm d}s \right|
=\left|\sum_{\sigma \in {\cal T}^h}  \sum_{i=1}^3 \int_{\partial_i \sigma} [\,(I-\mintpis)\,\eta^h]
\,\,\vv \,.\,\underline{\nu}_{\partial_i \sigma} \,{\rm d}s\right|
\nonumber \\
&=\left|\sum_{\sigma \in {\cal T}^h}
\sum_{i=1}^3 \int_{\partial_i \sigma} [\,(I-\mintpis)\,\eta^h]
\,\,[\,(\vI-\vP^h)\vv] \,.\,\underline{\nu}_{\partial_i \sigma} \,{\rm d}s \right|
\nonumber \\
&\leq\sum_{\sigma \in {\cal T}^h}  \sum_{i=1}^3 |\,(I-\mintpis)\,\eta^h|_{0,s,\partial_i \sigma}
\,\,|(\vI-\vP^h)\vv\,|_{0,s',\partial_i \sigma}
\nonumber \\
&\leq C \sum_{\sigma \in {\cal T}^h} h_\sigma \,|\del \eta^h|_{0,s,\sigma}
\,|\vv|_{1,s',\sigma}
\leq C\,h\,
|\del_h \eta^h|_{0,s,\Omega}
\,|\vv|_{1,s',\Omega}\,;
\label{T2bd1}
\end{align}
and hence the desired result (\ref{intbyparts}).

It immediately follows from (\ref{intbyparts}) that
\begin{align}
\left| \,(\eta^h, \del \cdot \vv) \right| \leq C\,|\del_h \eta^h|_{0,s,\Omega}
\,\|\vv\|_{1,s',\Omega}
\qquad \forall \
\eta^h \in N^h_0, \quad \forall \ \vv \in [W^{1,s'}_0(\Omega)]^2.
\label{NhPoinvetah}
\end{align}
Given any $\theta \in L^{s'}(\Omega)$, then there exists a $\vv \in [W^{1,s'}(\Omega)]^2$
such that
\begin{align}
\del\,.\,\vv = \theta \quad \mbox{a.e.\ in } \Omega, \qquad
\|\vv\|_{1,s',\Omega} \leq C\,|\theta|_{0,s',\Omega}.
\label{vvthetabd}
\end{align}
The result (\ref{vvthetabd}) is easily achieved by choosing $\vv = - \del z$, where
$- \Delta z = \theta'$ a.e.\ in $\Omega' \supset \Omega$ and $z=0$ on $\partial \Omega'$,
where $\theta'$ is the extension of $\theta$ from $\Omega$ to $\Omega'$
by zero and $\partial \Omega' \in C^{\infty}$.
Combining (\ref{NhPoinvetah}) and (\ref{vvthetabd}) yields the desired result
(\ref{NhPoin}).
\end{proof}

\begin{lmm}\label{Nhcompact}
Given $\{\eta^h\}_{h>0}$, with $\eta^h \in N^h_0$, such that for an $s \in (4,\infty)$
\begin{align}
|\del_h \eta^h|_{0,s,\Omega} \leq C\,;
\label{deletahstab}
\end{align}
then there exists
a subsequence of $\{\eta^h\}_{h>0}$, (not indicated),
and
an $\eta \in W^{1,s}_0(\Omega)$ such that
as $h \rightarrow 0$
\eqlabstart
\begin{alignat}{2}
\del_h \eta^h &\rightarrow \del \eta \qquad &&\mbox{weakly in } [L^s(\Omega)]^2,
\label{weakdeletah}\\
\eta^h &\rightarrow \eta \qquad &&\mbox{strongly in } L^s(\Omega),
\label{strongetah}\\
\del_h \eta^h &\rightarrow \del \eta \qquad &&\mbox{strongly in }
[\,[H^{\frac{1}{2}}(\Omega)]^\star]^2.
\label{strongdeletah}
\end{alignat}
\eqlabend
\end{lmm}
\begin{proof}
It follows immediately from (\ref{deletahstab}) and (\ref{NhPoin}) that there exist
an $\eta \in L^s(\Omega)$ and a
$\vd \in [L^s(\Omega)]^2$, and a subsequence of $\{\eta^h\}_{h>0}$ (not indicated)
such that as $h \rightarrow 0$
\begin{align}
&\eta^h \rightarrow \eta \qquad \mbox{weakly in } L^s(\Omega), \qquad \qquad
\del_h \eta^h \rightarrow \vd \qquad \mbox{weakly in } [L^s(\Omega)]^2.
\label{weakdeletahd}
\end{align}
Passing to the limit $h \rightarrow 0$ in (\ref{intbyparts}) for the subsequence
we deduce that
\begin{align}
\eta \in W^{1,s}_0(\Omega) \qquad \mbox{and} \qquad \vd= \del \eta\,.
\label{vdident}
\end{align}
Hence the desired result (\ref{weakdeletah}) follows from combining
(\ref{weakdeletahd}) and (\ref{vdident}).

We now introduce $\hat \eta^h \in U^h_0$ such that
\begin{align}
(\del \hat \eta^h, \del \chi^h)
= (\del_h \eta^h, \del \chi^h) \qquad \forall \ \chi^h \in U^h_0\,.
\label{hatetah}
\end{align}
It follows from (\ref{Poin}), (\ref{hatetah}) and (\ref{deletahstab}) that
\begin{align}
\| \hat \eta^h\|_{1,\Omega} \leq C\,|\del \hat \eta^h|_{0,\Omega}
\leq C\,|\del_h \eta^h|_{0,\Omega} \leq C\,.
\label{hatetahH1}
\end{align}
We deduce from (\ref{hatetahH1}) that there exists
a further subsequence of $\{\eta^h \}_{h>0}$
(not indicated) such that as $h \rightarrow 0$
\eqlabstart
\begin{alignat}{3}
\del \hat \eta^h
& \rightarrow \del \hat \eta
\qquad &&
\mbox{weakly in } [L^2(\Omega)]^2,
\label{wwcondelhatetah}\\
\hat \eta^h&\rightarrow \hat \eta  \qquad &&
\mbox{strongly in } L^\kappa(\Omega), \qquad \forall \ \kappa \in [1,\infty);
\label{wsconhatetah}
\end{alignat}
\eqlabend
where $\hat \eta \in H^1_0(\Omega)$.

As $\Omega$ is polygonal,
it follows from Grisvard \cite[Chapter 4]{Grisvard} that given $\theta \in L^{s'}(\Omega)$
for some $s' \in (1, \frac{4}{3})$, then
there exists a unique $z \in W^{2,s'}(\Omega)$ such that
\begin{align}
&- \Delta z = \theta \quad \mbox{a.e.\ in } \Omega, \qquad z=0 \quad \mbox{on } \partial \Omega;
\qquad\mbox{and}\qquad
\|z\|_{2,s',\Omega} \leq C \,|\theta|_{0,s',\Omega}\,.
\label{zthetamubd}
\end{align}
It follows from (\ref{zthetamubd}) and (\ref{hatetah}) that
\begin{align}
&( \hat \eta^h - \eta^h, \theta) =
[ \,(\del \hat \eta^h  - \del_h \eta^h, \del[(I-\pi^h_U)z]\,)]
+ [ \,( \del_h \eta^h, \del z) + (\eta^h, \Delta z)\,]
=: T_1 + T_2\,.
\label{T1T2}
\end{align}
We deduce from (\ref{hatetahH1}), (\ref{piconv0}) and (\ref{zthetamubd}) that
\eqlabstart
\begin{align}
|T_1| \leq C\,|(I-\pi^h_U)z|_{1,\Omega}
\leq C\,h^{2(1-\frac{1}{s'})}\,|z|_{2,s',\Omega}
\leq C\,h^{\frac{2}{s}}\,|\theta|_{0,s',\Omega}\,,
\label{T1bd}
\end{align}
and from (\ref{intbyparts}) with $\vv = - \del z$, (\ref{deletahstab}) and (\ref{zthetamubd}) that
\begin{align}
|T_2| \leq C \,h\,|\del_h \eta^h|_{0,s,\Omega}\,\|z\|_{2,s',\Omega}
\leq C\,h\,|\theta|_{0,s',\Omega}\,.
\label{T2bd}
\end{align}
\eqlabend
Hence combining (\ref{T1T2}) and (\ref{T1bd},b) yields
that
\begin{align}
|\,(\hat \eta^h - \eta^h, \theta)| \leq C\,h^{\frac{2}{s}}\,
|\theta|_{0,s',\Omega}
\qquad \forall \ \theta \in L^{s'}\!(\Omega)\,.
\label{T1T2bd}
\end{align}
It follows immediately from (\ref{T1T2bd}) that
\begin{align}
|\hat \eta^h - \eta^h|_{0,s,\Omega} \rightarrow 0 \qquad \mbox{as } \quad
h \rightarrow 0\,.
\label{T1T2bd1}
\end{align}
The desired result (\ref{strongetah}) then follows from
(\ref{wsconhatetah}) and (\ref{T1T2bd1}), as (\ref{weakdeletahd}) implies that $\hat \eta = \eta$.

Finally, we need to prove (\ref{strongdeletah}). First, we note that
\begin{align}
(\del \eta - \del_h \eta^h, \vv) = - (\eta-\eta^h, \del\,.\,\vv) -
\left[ (\eta^h, \del\,.\,\vv) + (\del_h \eta^h, \vv) \right] \qquad \forall \ \vv \in
[H^1(\Omega)]^2.
\label{intbypartsa}
\end{align}
Hence it follows from   (\ref{intbypartsa}) and (\ref{intbyparts}) that
\begin{align}
\left| (\del \eta - \del_h \eta^h, \vv) \right|
\leq C\,\left[ |\eta-\eta^h|_{0,\Omega} + h\,|\del_h \eta^h|_{0,\Omega} \right]
|\vv|_{1,\Omega} \qquad \forall \ \vv \in [H^1(\Omega)]^2.
\label{intbypartsb}
\end{align}
Therefore (\ref{intbypartsb}), (\ref{strongetah}) and (\ref{deletahstab})
yield for the subsequence that
\begin{alignat}{2}
\del_h \eta^h &\rightarrow \del \eta \qquad &&\mbox{strongly in }
[\,[H^{1}(\Omega)]^\star]^2 \mbox{ as $h \rightarrow 0$}.
\label{strongdeletah1star}
\end{alignat}
The desired result (\ref{strongdeletah}) follows immediately
from (\ref{Hhalfstarnorm}), (\ref{strongdeletah1star}) and (\ref{deletahstab}).
\end{proof}

\section{Convergence}
\label{conv}
\setcounter{equation}{0}

\subsection{Convergence of (Q$^{h,\tau}_{r}$) to (Q$^{\tau}_r$)}
\label{secconvB}

Similarly to (\ref{calMh}), 
we introduce for all $\eta \in C(\overline{\Omega})$
\begin{align}
{\mathfrak M}^{n} (\eta)(\vx) \: := \left\{\begin{array}{ll}
M_\eps(\eta)(\vx) &\qquad \mbox{in case (i),}\\[2mm]
k(\vx)\,\widehat M(\eta(\vx)+b_e(t_n)) &\qquad \mbox{in case (ii),}\\[2mm]
k(\vx)&\qquad \mbox{in case (iii)}
\end{array}
\right.
\qquad \mbox{for a.e.\ } \vx \in \Omega.
\label{calM}
\end{align}
We note from (\ref{calM}), (\ref{Meps}) and Assumption (A1) that
there exist ${\mathfrak M}_{\rm min},\,{\mathfrak M}_{\rm max} \in {\mathbb R}$ such that
for $n=1,\ldots,N$
\begin{align}
0 < {\mathfrak M}_{\rm min} \leq {\mathfrak M}^{n}(\eta)(\vx) \leq {\mathfrak M}_{\rm max}
\qquad \forall \ \eta \in C(\overline{\Omega}), \qquad
\mbox{for a.e.\ } \vx \in \Omega\,.
\label{calMb}
\end{align}

For the purposes of the convergence analysis in this subsection, we introduce
for a given $r>1$:

\noindent
{\bf (Q$^{\tau}_r$)}
For $n = 1,\ldots, N$, find $w^n_r \in W^{1,p}_0(\Omega)$ and
$\vq^n_r \in [L^{r}(\Omega)]^2$ such that
\eqlabstart
\begin{alignat}{2}
{\cal A}\left(\frac{w^n_r-w^{n-1}_r}{\tau_n},\eta\right) - (\vq^n_r,\del\eta) &= ({\cal F}^n
,\eta)
\quad \;&&\forall \ \eta \in W^{1,p}_0(\Omega),\label{Qt1r}\\
( {\mathfrak M}^n(w^n_r)\,|\vq^n_r|^{r-2}\vq^n_r,\vv) + (\del w^n_r,\vv)  &=0
\quad \;
&&\forall \ \vv\in [L^{r}(\Omega)]^2\,;
\label{Qt2r}
\end{alignat}
\eqlabend
where $w^0_r = w^0$.

\begin{thrm}\label{hconthmB}
Let the Assumptions (A1) and(A2) hold.
For any fixed $r \in (1,\tfrac{4}{3})$ and fixed time partition $\{\tau_n\}_{n=1}^N$,
and for all regular partitionings ${\cal T}^h$ of $\Omega$,
there exists a subsequence of $\{\{W^n_{r},$ $\vQ^n_{r}\}_{n=1}^N\}_{h>0}$ (not indicated),
where $\{W^n_{r},\vQ^n_{r}\}_{n=1}^N$ solves (Q$^{h,\tau}_{r}$), (\ref{Qthr1B},b),
such that as $h \rightarrow 0$, for any $s \in [1,\infty)$,
\eqlabstart
\begin{alignat}{3}
W^{n}_{r}, \, P^h W^n_r&\rightarrow w^n_r  \qquad &&
\mbox{strongly in } L^p(\Omega),
\qquad && n=1,\ldots, N,
\label{wsconBr}
\\
{\mathfrak M}^{h,n}(P^h W^{n}_{r}) &\rightarrow {\mathfrak M}^n(w^n_r)  \qquad &&
\mbox{strongly in } L^s(\Omega),
\qquad && n=1,\ldots, N,
\label{MsconBr}
\\
\del_h W^n_r &\rightarrow \del w^n_r \qquad
&&\mbox{weakly in } [L^p(\Omega)]^2,
\qquad && n=1,\ldots, N,
\label{wdelwconBr}\\
\del_h W^n_r &\rightarrow \del w^n_r \qquad &&\mbox{strongly in }
[\,[H^{\frac{1}{2}}(\Omega)]^\star]^2,\qquad && n=1,\ldots, N,
\label{strongdelWhalfstar} \\
\vQ^n_{r} &\rightarrow \vq^n_r \qquad &&
\mbox{weakly in } [L^r(\Omega)]^2,
\qquad && n=1,\ldots, N;
\label{qwconBr}
\end{alignat}
\eqlabend
where $\{w^n_r,\vq^n_r\}_{n=1}^N$ is a solution of (Q$^\tau_r$), (\ref{Qt1r},b).

In addition, we have that
\begin{align}
\max_{n=0,\ldots, N}
\|w^n_{r}\|_{\cal A} + \sum_{n=1}^N \|w^n_{r}-w^{n-1}_{r}\|_{\cal A}^2 +
\sum_{n=1}^N \tau_n\,|\vq^n_{r}|_{0,r,\Omega}^r
+\left( \sum_{n=1}^N \tau_n |\del w^n_{r}|_{0,p,\Omega}^p\right)^{\frac{1}{p}}
&\leq C.
\label{energy4}
\end{align}
Moreover, in case (iii) the solution of (Q$^\tau_r$) is unique, and so the whole sequence
converges in (\ref{wsconBr}--e).
\end{thrm}
\begin{proof}
The desired subsequence weak convergence result (\ref{qwconBr}) follows immediately from the
bound on $\{\vQ^n_{r}\}_{n=1}^N$ in (\ref{energy3}), on noting that the time
partition $\{\tau_n\}_{n=1}^N$ is fixed.
It follows from 
(\ref{energy3}) that
\begin{align}
|\del_h W^n_{r}|_{0,p,\Omega} \leq C(\tau_n^{-1}), \qquad n=1,\ldots, N\,.
\label{delhWnBrL2}
\end{align}
The desired results (\ref{wsconBr},c,d) then follow immediately from
(\ref{delhWnBrL2}), Lemma \ref{Nhcompact}
and (\ref{Phstab},b)
on extracting a further subsequence (not indicated).
On noting
that ${\mathfrak M}^n(\cdot)$ is well-defined on $S^h$
and is continuous with respect to its argument,
it follows from (\ref{wsconBr})
for a further subsequence of $\{\{W^n_{r}\}_{n=1}^N\}_{h>0}$ (not indicated)
that as $h \rightarrow 0$, for $n=1,\ldots, N$,
\begin{align}
P^h W^n_{r} \rightarrow w^n_r \quad \mbox{a.e.\ in } \Omega
\quad \Rightarrow \quad {\mathfrak M}^n(P^h W^n_{r}) \rightarrow {\mathfrak M}^n(w^n_r)
\quad \mbox{a.e.\ in } \Omega\,.
\label{MWwae}
\end{align}
It follows from (\ref{MWwae}), 
(\ref{calMb}) and Lebesgue's general
convergence theorem that as $h \rightarrow 0$
for any $s \in [1,\infty)$
\begin{align}
{\mathfrak M}^n(P^h W^n_{r}) \rightarrow {\mathfrak M}^n(w^n_r) \qquad \mbox{strongly in } L^s(\Omega),
\qquad n=1,\ldots, N\,.
\label{MsconBrs}
\end{align}
Combining (\ref{calM}), (\ref{calMh}), (\ref{Mepshapprox}), (\ref{Pconv}),
(\ref{piconv})
and (\ref{MsconBrs}) yields the desired result (\ref{MsconBr}).

We now need to establish that $\{w^n_r,\vq^n_r\}_{n=1}^N$ solve (Q$^\tau_r$), (\ref{Qt1r},b).
For any $\eta \in C^\infty_0(\Omega)$, we choose  $\eta^h = \pi^h_N \eta$ in (\ref{Qthr1B})
and now pass to the limit $h \rightarrow 0$ for the subsequence,
on noting (\ref{calAh}), (\ref{calA}), (\ref{bilformc},b), (\ref{wsconBr},d,e)
and (\ref{piconv}),
to obtain (\ref{Qt1r}) for all $\eta \in C^\infty_0(\Omega)$.
Noting that $C^\infty_0(\Omega)$ is dense in $W^{1,p}_0(\Omega)$, (\ref{calA}), (\ref{bilformc},b)
and that $w^n_r \in W^{1,p}_0(\Omega)$,
$\vq^n_r \in [L^r(\Omega)]^d$ and ${\cal F}^n \in L^2(\Omega)$, $n=1,\ldots,N$,
yields the desired result (\ref{Qt1r}).

For any $\vv \in [C^\infty(\overline{\Omega})]^2$, we choose  $\vv^h = \vQ^n_{r}-\vP^h \vv$
in (\ref{Qthr2B}), 
and then try to pass to the limit for the subsequence as $h \rightarrow 0$.
First, we note from (\ref{Qthr1B}) with $\eta^h = W^n_r$ and
(\ref{aederiv}) 
that
for $n=1,\ldots, N$
\begin{align}
(\del_h W^n_{r},\vP^h \vv) &=
(\del_h W^n_{r},\vQ^n_r) +
({\mathfrak M}^{h,n}(P^h W^n_{r})\,|\vQ^n_{r}|^{r-2}\,\vQ^n_{r},\vQ^n_{r}-\vP^h \vv)
\nonumber \\
&\geq
{\cal A}^h( \frac{W^n_r-W^{n-1}_r}{\tau_n},W^n_r) - ({\cal F}^n, W^n_r)
+
({\mathfrak M}^{h,n}(P^h W^n_{r}),|\vP^h \vv|^{r-2} \vP^h \vv,
\vQ^n_{r}-\vP^h \vv)
\,.
\label{Qthr2BtoQtr2b}
\end{align}
Passing to the limit $h \rightarrow 0$ for the subsequence in (\ref{Qthr2BtoQtr2b})
yields,
on noting (\ref{wsconBr}--e), (\ref{Pconv}), (\ref{calAh}), (\ref{calA})
and (\ref{bilformc},b), 
for $n=1,\ldots, N$
that
\begin{align}
(\del w^n_{r},\vv)&\geq
{\cal A}(\frac{w^n-w^{n-1}_r}{\tau_n},w^n_r)
- ({\cal F}^n, w^n_r)
+
({\mathfrak M}^{n}(w^n_{r}),|\vv|^{r-2} \vv,
\vq^n_{r}-\vv)
\,.
\label{Qtrcomb}
\end{align}
It follows from (\ref{Qtrcomb}) and (\ref{Qt1r}) with $\eta=w^n_r$ that
for $n=1,\ldots, N$
\begin{align}
(\del w^n_r, \vv-\vq^n_r) \geq ({\mathfrak M}^n(w^n_r) \,|\vv|^{r-2}\,\vv, \vq^n_r - \vv)
\qquad \forall \ \vv \in [C^\infty(\overline{\Omega})]^2\,.
\label{Qthr2BtoQtr2d}
\end{align}
As $w^n_r \in C(\overline{\Omega})$, ${\mathfrak M}^n(w^n_r) \in L^{\infty}(\Omega)$
and $\vq^n_r \in [L^r(\Omega)]^2$,
it follows 
that
(\ref{Qthr2BtoQtr2d}) holds true for all $\vv \in [L^r(\Omega)]^2$.
For any fixed $\vz \in [L^r(\Omega)]^2$, choosing $\vv = \vq^n_r \pm \alpha \vz$ with
$\alpha \in {\mathbb R}_{>0}$ in (\ref{Qthr2BtoQtr2d}) and letting $\alpha \rightarrow 0$
yields the desired result (\ref{Qt2r}) on repeating the above for any
$\vz \in [L^r(\Omega)]^2$.

In addition, it follows from $W^0_r = \pi^h_N w^0$ and (\ref{piconv}) that $w^0_r = w^0$.
Therefore $\{w^n_r,\vq^n_r\}_{n=1}^N$ is a solution of (Q$^{\tau}_r$),
(\ref{Qt1r},b). It follows from (\ref{energy3}), (\ref{calAh}), (\ref{bilformc},b), (\ref{calA}),
(\ref{wsconBr},c,d,e) and (\ref{aederiv}) that (\ref{energy4}) holds.

Finally, it is a simple matter to establish the uniqueness of the solution of (Q$^\tau_r$)
in case (iii).
\end{proof}

\begin{crllr}
\label{hatWnrcorr}
Let the Assumptions of Theorem \ref{hconthmB} hold.
For $n=1,\ldots,N$ let $\hat W^n_r \in U^h_0$ be such that
\begin{align}
(\del \hat W^n_r, \del \chi^h) = (\del_h W^n_r, \del \chi^h) \qquad
\forall \ \chi^h \in U^h_0.
\label{hatWnr}
\end{align}
Then
there exists a further subsequence of $\{\{W^n_{r},$ $\vQ^n_{r}\}_{n=1}^N\}_{h>0}$ (not indicated),
where $\{W^n_{r},\vQ^n_{r}\}_{n=1}^N$ solves (Q$^{h,\tau}_{r}$), (\ref{Qthr1B},b),
such that as $h \rightarrow 0$
\eqlabstart
\begin{alignat}{3}
\hat W^{n}_{r} &\rightarrow w^n_r  \qquad &&
\mbox{strongly in } L^s(\Omega),
\qquad && n=1,\ldots, N, \qquad \forall \ s \in [1,\infty),
\label{wsconhatWnr}
\\
\del \hat  W^n_r &\rightarrow \del w^n_r \qquad
&&\mbox{weakly in } [L^2(\Omega)]^2,
\qquad && n=1,\ldots, N;
\label{wdelwconhatWnr}
\end{alignat}
\eqlabend
where $\{w^n_r,\vq^n_r\}_{n=1}^N$ is a solution of (Q$^\tau_r$), (\ref{Qt1r},b).
In case (iii) the whole sequence converges in (\ref{wsconhatWnr},b)
as the solution of (Q$^\tau_r$) is unique.
\end{crllr}
\begin{proof}
The proof follows immediately from (\ref{delhWnBrL2}), (\ref{wsconBr},c), (\ref{hatWnr}),
(\ref{deletahstab}), (\ref{weakdeletah},b), (\ref{hatetah}) and (\ref{wwcondelhatetah},b)
on noting that $\hat \eta = \eta$.

\end{proof}

\subsection{Convergence of (Q$^\tau_r$) to (Q) in case (i)}

It follows from (\ref{Qt1r}), (\ref{energy4}), (\ref{calA}) and (\ref{calFn})
in the growing sandpile case that for $n=1,\ldots,N$
\begin{align}
\tau_n\, | (\vq^n_r, \del \eta) | \leq C\,|\eta|_{0,\Omega}
\qquad \forall \ \eta \in C^\infty_0(\Omega).
\label{divqnrbd}
\end{align}
Hence, for a fixed time partition $\{\tau_n\}_{n=1}^N$,
the distributional divergence of $\vq^n_r$ belongs $L^2(\Omega)$, $n=1,\ldots,N$.
Therefore,
on recalling (\ref{vVs}),
(Q$^\tau_r$), (\ref{Qt1r},b), can be reformulated for a given $r \in (1, \frac{4}{3})$ as:
\vspace{2mm}

\noindent
{\bf (Q$^{\tau}_r$)}
For $n = 1,\ldots, N$, find $w^n_r \in W^{1,p}_0(\Omega)$ and
$\vq^n_r \in \vV^r(\Omega)$ such that
\eqlabstart
\begin{alignat}{2}
\left(\frac{w^n_r-w^{n-1}_r}{\tau_n},\eta\right) + (\del\,.\,\vq^n_r,\eta) &=
(f^n
,\eta)
\quad \;&&\forall \ \eta \in L^2(\Omega),\label{Qt1rsand}\\
(M_\eps(w^n_r)\,|\vq^n_r|^{r-2}\vq^n_r,\vv) - (w^n_r,\del\,.\,\vv)  &=0
\quad \;
&&\forall \ \vv\in \vV^r(\Omega)\,;
\label{Qt2rsand}
\end{alignat}
\eqlabend
where $w^0_r = w_0^\epsilon$.

The above is the formulation of (Q$^\tau_r$) in Barrett and Prigozhin \cite[(3.24a,b)]{sandqvi}.
On recalling (\ref{vVM}),
we state the discrete time approximation of the mixed formulation
of the growing sandpile problem; that is, the $r \rightarrow 1$ limit of (Q$^\tau_r$):
\vspace{2mm}

\noindent
{\bf (Q$^{\tau}$)}
For $n = 1,\ldots, N$, find $w^n \in W^{1,\infty}_0(\Omega)$ and
$\vq^n \in \vV^{\cal M}(\Omega)$ such that
\eqlabstart
\begin{alignat}{2}
\left(\frac{w^n-w^{n-1}}{\tau_n},\eta\right) + (\del\,.\,\vq^n,\eta) &= (f^n
,\eta)
\quad \;&&\forall \ \eta \in L^2(\Omega),\label{Qt1}\\
\langle |\vv|-|\vq^n|,M_\epsilon(w^n)\rangle_{C(\overline{\Omega})} - (\del\,.\,(\vv-\vq^n),
w^n)  &\geq0
\quad \;
&&\forall \ \vv\in \vV^{\cal M}(\Omega)\,;
\label{Qt2}
\end{alignat}
\eqlabend
where $w^0 = w_0^\epsilon$.

Similarly to (\ref{Kset}), we introduce for
$\chi \in W^{1,\infty}_0(\Omega)$
the closed convex non-empty set
\begin{align}
K_\epsilon(\chi) := \{ \eta \in W^{1,\infty}_0(\Omega) : |\del \eta| \leq M_\epsilon(\chi)
\quad \mbox{a.e.\ on } \Omega\}\,.
\label{K}
\end{align}
Then associated with (Q$^{\tau}$) is the corresponding approximation of the primal
quasi-variational inequality:

\noindent
{\bf (P$^{\tau}$)}
For $n = 1,\ldots, N$, find $w^n \in K_\epsilon(w^n)$ such that
\begin{align}
&\left(\frac{w^n-w^{n-1}}{\tau_n},\eta-w^n\right)
\geq (f^n ,\eta-w^n)
\qquad \forall \ \eta \in K_\epsilon(w^n)\,,\label{Pt}
\end{align}
where $w^0 = w_0^\epsilon$.

\vspace{2mm}

Similarly to \cite{sandqvi}, for our convergence results we require extra assumptions. 

\vspace{2mm}

\noindent
{\bf (A3)} $\Omega$ is a strictly star-shaped domain.

\vspace{2mm}

\noindent
{\bf (A4)}
$w^\epsilon_0 \geq 0$ and
$f \in L^\infty(0,T;L^2(\Omega))$.

\vspace{2mm}

\begin{thrm}\label{hconthmr}
Let the Assumptions (A1), (A2) and (A3) hold.
For any fixed time partition $\{\tau_n\}_{n=1}^N$, 
there exists a subsequence of $\{\{w^n_{r},\vq^n_{r}\}_{n=1}^N \}_{r\in(1,\frac{4}{3})}$ (not indicated),
where $\{w^n_{r},\vq^n_{r}\}_{n=1}^N$ solves (Q$^{\tau}_{r}$), (\ref{Qt1rsand},b),
such that as $r \rightarrow 1$ 
\eqlabstart
\begin{alignat}{3}
w^{n}_{r}&\rightarrow w^n  \qquad &&
\mbox{strongly in } C(\overline{\Omega}),
\qquad && n=0,\ldots, N,
\label{wsconr}
\\
M_\epsilon( w^{n}_{r}) &\rightarrow M_\epsilon(w^n)  \qquad &&
\mbox{strongly in } C(\overline{\Omega}),
\qquad && n=0,\ldots, N,
\label{Msconr}
\\
\vq^n_{r} &\rightarrow \vq^n \qquad &&
\mbox{weakly in } [{\cal M}(\overline{\Omega})]^2,
\qquad && n=1,\ldots, N,
\label{qwconr}\\
\del\,.\,\vq_{r}^n &\rightarrow \del\,.\,\vq^n \qquad &&
\mbox{weakly in } L^2(\Omega),
\qquad && n=1,\ldots, N;
\label{hdelqconr}
\end{alignat}
\eqlabend
where $\{w^n,\vq^n\}_{n=1}^N$ is a solution of (Q$^\tau$), (\ref{Qt1},b).
\end{thrm}
\begin{proof}
See the proof of Theorem 3.4 in Barrett and Prigozhin \cite{sandqvi}.
We note that the convexity of $\Omega$ and the restriction of $\tau \in (0,\frac{1}{2}]$
were also assumed there,
as these were required solely to establish the existence of a solution $\{w^n_r, \vq^n_r\}_{n=1}^N$
to (Q$^{\tau}_r$), see Theorem 3.3 in Barrett and Prigozhin \cite{sandqvi}.
These constraints on $\Omega$ and $\tau$ are not required here, see Theorem \ref{hconthmB} above.
In addition, as the time partition $\{\tau_n\}_{n=1}^N$ is fixed, the bound
$|\del\,.\,\vq^r|_{0,\Omega} \leq C(\tau_n^{-1})$, $n=1,\ldots,N$,
which immediately follows from (\ref{divqnrbd}) is adequate to establish (\ref{hdelqconr}).
Therefore the bound on $\del\,.\,\vq^n_r$ in Barrett and Prigozhin \cite[(3.47)]{sandqvi} is not necessary.
\end{proof}

Next, we note the following result.

\begin{thrm}
\label{QPtaueq}
Let the Assumptions (A1), (A2) and (A3) hold.
If $\{w^n,\vq^n\}_{n=1}^N$ is a solution of (Q$^{\tau}$), (\ref{Qt1},b), then
$\{w^n\}_{n=1}^N$ solves (P$^{\tau}$), (\ref{Pt}), and
\begin{align}
w^n \geq w^{n-1} \qquad n=1,\ldots, N\,.
\label{wnmono}
\end{align}
\end{thrm}
\begin{proof}
See the proof of Theorem 3.6 in Barrett and Prigozhin \cite{sandqvi}.
\end{proof}

We introduce the following notation for $t \in (t_{n-1},t_n]$, $n = 1,\ldots, N$,
\begin{align}
f^{\tau,+}(\cdot,t) &:= f^n(\cdot)\, \qquad
w^{\tau}(\cdot,t) := \frac{(t-t_{n-1})}{\tau_n}\,w^n(\cdot)
+\frac{(t_n-t)}{\tau_n}\,w^{n-1}(\cdot)\,,
\nonumber \\
w^{\tau,+}(\cdot,t) &:= w^n(\cdot),
\qquad
w^{\tau,-}(\cdot,t) := w^{n-1}(\cdot),
\qquad
\vq^{\tau,+}(\cdot,t) := \vq^n(\cdot)\,.
 \label{wqtau+}
\end{align}

We now introduce the weak mixed formulation of the growing sandpile problem:

\noindent
{\bf (Q)}
Find $w \in L^\infty(0,T;W^{1,\infty}_0(\Omega))\cap W^{1,\infty}(0,T;
[C^1_0(\overline{\Omega})]^\star)$ and
$\vq \in L^\infty(0,T;$ $[{\cal M}(\overline{\Omega})]^2)$ such that
\eqlabstart
\begin{align}
&\int_0^T \left[
\langle\frac{\partial w}{\partial t},\eta \rangle_{C^1_0(\overline{\Omega})} -
\langle \vq,\del \eta\rangle_{C(\overline{\Omega})} -(f
,\eta) \right] \,{\rm d}t =0
\qquad
\forall \ \eta \in L^1(0,T;C^1_0(\overline{\Omega})),\label{Q1}\\
&\int_0^T \left[ \langle |\vv|-|\vq|,M_\epsilon(w)\rangle_{C(\overline{\Omega})} -
(\del\,.\,\vv-f,
w)\right] \,{\rm d}t  
\geq \tfrac{1}{2} \left[\, |w(\cdot,T)|^2_{0,\Omega} - |w^\epsilon_0(\cdot)|^2_{0,\Omega}
\, \right]
\qquad \forall \ \vv\in L^1(0,T;\vV^{\cal M}(\Omega))\,;
\label{Q2}
\end{align}
\eqlabend
where $w(\cdot,0) = w_0^\epsilon(\cdot)$.

Associated with (Q) is the corresponding primal
quasi-variational inequality:

\noindent
{\bf (P)}
Find $w \in L^\infty(0,T;K_\epsilon(w))\cap
W^{1,\infty}(0,T;[C^1_0(\overline{\Omega})]^\star)$ such that
\begin{align}
&\int_0^T \left[
\langle\frac{\partial w}{\partial t},\eta \rangle_{C^1_0(\overline{\Omega})} -
(f ,\eta-w) \right] \,{\rm d}t
\geq
\tfrac{1}{2} \left[ \,|w(\cdot,T)|^2_{0,\Omega} - |w^\epsilon_0(\cdot)|^2_{0,\Omega}
\,\right]
\qquad
\forall \ \eta \in L^1(0,T;K_\epsilon(w)\cap C^1_0(\overline{\Omega})),\label{P}
\end{align}
where $w(\cdot,0) = w_0^\epsilon(\cdot)$.

For the reasoning behind the formulations (Q) and (P), and the Assumption (A4); see Remarks 3.1
and 3.9 in Barrett and Prigozhin \cite{sandqvi}.

\begin{thrm}\label{conthm}
Let the Assumptions (A1), (A2), (A3) and (A4) hold.
For all time partitions $\{\tau_n\}_{n=1}^N$,
there exists a subsequence of $\{\{w^n,\vq^n\}_{n=1}^N\}_{\tau>0}$ (not indicated),
where $\{w^n,\vq^n\}_{n=1}^N$ solves (Q$^{\tau}$), (\ref{Qt1},b),
such that as $\tau \rightarrow 0$
\eqlabstart
\begin{alignat}{2}
w^{\tau},\,w^{\tau,\pm}&\rightarrow w  \qquad &&
\mbox{weak$^\star$ in } L^\infty(0,T;W^{1,\infty}(\Omega)),
\label{wgradwcont}\\
\frac{\partial w^{\tau}}{\partial t}&\rightarrow \frac{\partial w}{\partial t}  \qquad &&
\mbox{weak$^\star$ in } L^\infty(0,T;[C^{1}_0(\overline{\Omega})]^\star),
\label{wtwcont}\\
w^{\tau} &\rightarrow w  \qquad &&
\mbox{strongly in } C([0,T];C(\overline{\Omega})),
\label{wscont}
\\
w^{\tau,\pm}&\rightarrow w  \qquad &&
\mbox{strongly in } L^2(0,T;C(\overline{\Omega})),
\label{wsL2}
\\
M_\epsilon(w^{\tau})&\rightarrow M_\epsilon(w)  \qquad &&
\mbox{strongly in } C([0,T];C(\overline{\Omega})),
\label{Mwscont}
\\
M_\epsilon(w^{\tau,\pm}) &\rightarrow M_\epsilon(w)  \qquad &&
\mbox{strongly in } L^2(0,T;C(\overline{\Omega})),
\label{MwsL2}
\\
\vq^{\tau,+} &\rightarrow \vq \qquad &&
\mbox{weak$^\star$ in } L^\infty(0,T;[{\cal M}(\overline{\Omega})]^d);
\label{qwcont}
\end{alignat}
\eqlabend
where $\{w,\vq\}$ is a solution of (Q), (\ref{Q1},b).
Moreover, $w$ solves (P), (\ref{P}).
\end{thrm}
\begin{proof}
See the proof of Theorem 3.8 in Barrett and Prigozhin \cite{sandqvi}.
\end{proof}

\subsection{Convergence of (Q$^{\tau}_r$) to (Q) in case (ii)}

In the cylindrical superconductor case, on noting (\ref{calM}), (\ref{Qt2r}) becomes
\begin{align}
( k\, \widehat M(w^n_r+b_e(t_n))\,|\vq^n_r|^{r-2}\vq^n_r,\vv) + (\del w^n_r,\vv)  &=0
\qquad \forall \ \vv\in [L^{r}(\Omega)]^2\,,
\label{Qt2rii}
\end{align}
which can be rewritten, on noting (\ref{jcsc}), as
\begin{align}
( k\,|\vq^n_r|^{r-2}\vq^n_r,\vv) +
([\widehat M(w^n_r+b_e(t_n))]^{-1}\del w^n_r,\vv)  &=0
\qquad \forall \ \vv\in [L^{r}(\Omega)]^2\,.
\label{Qt2riiM}
\end{align}
With $\widehat F \in C({\mathbb R};{\mathbb R})$ such that
\begin{align}
(\hat M_0)^{-1} \geq
\widehat F'(s) = [\widehat M (s)]^{-1} \geq (\hat M_1)^{-1} >0\qquad \mbox{and} \qquad \widehat F(0)=0,
\label{wF}
\end{align}
where we have noted (\ref{jcsc}),
(\ref{Qt2riiM}) can be rewritten as
\begin{align}
( k\,|\vq^n_r|^{r-2}\vq^n_r,\vv) +
(\del [\widehat F(w^n_r+b_e(t_n)) - \widehat F(b_e(t_n))] ,\vv)  &=0
\qquad \forall \ \vv\in [L^{r}(\Omega)]^2\,.
\label{Qt2riiF}
\end{align}
Then similarly to (\ref{Qt1rsand},b), on noting the analogue of (\ref{divqnrbd}),
(\ref{fn}), (\ref{calFn}) and (\ref{Qt2riiF}),
(Q$^{\tau}_r$), (\ref{Qt1r},b), in the cylindrical superconductor case can be reformulated
for a given $r \in (1,\frac{4}{3})$ as:

\noindent
{\bf (Q$^{\tau}_r$)}
For $n = 1,\ldots, N$, find $w^n_r \in W^{1,p}_0(\Omega)$ and
$\vq^n_r \in \vV^r(\Omega)$ such that
\eqlabstart
\begin{alignat}{2}
\left(\frac{(w^n_r+b_e(t_n))-(w^{n-1}_r+b_e(t_{n-1}))}
{\tau_n},\eta\right) + (\del\,.\,\vq^n_r,\eta)
&=  0
\qquad&&\forall \ \eta \in L^2(\Omega),\label{Qt1rsc}\\
( k\,|\vq^n_r|^{r-2}\vq^n_r,\vv)
-(\widehat F(w^n_r+b_e(t_n)) - \widehat F(b_e(t_n)),\del\,.\,\vv)  &=0
\qquad &&\forall \ \vv\in \vV^{r}(\Omega)\,.
\label{Qt2rsc}
\end{alignat}
\eqlabend

It is now a simple matter to establish the uniqueness of $\{w^n_r,\vq^n_r\}_{n=1}^N$
solving (Q$^{\tau}_r$), (\ref{Qt1rsc},b), by exploiting (\ref{aederiv})
and the monotonicity of $\widehat{F}$, recall (\ref{wF}).
In addition, we have the following stability result.

\begin{lmm}\label{hconthmBsc}
Let the Assumptions (A1) and(A2) hold.
For any fixed $r \in (1,\tfrac{4}{3})$ and time partition $\{\tau_n\}_{n=1}^N$,
the unique solution $\{w^n_r,\vq^n_r\}_{n=1}^N$ of (Q$^{\tau}_r$), (\ref{Qt1rsc},b),
in addition to satisfying
(\ref{energy4}) with $\|\cdot\|_{\cal A} \equiv |\cdot|_{0,\Omega}$ satisfies
\begin{align}
\textstyle \frac{r-1}{r}\,\displaystyle \max_{n=1,\dots, N}
(k,|\vq^n_r|^r) +
\sum_{n=1}^N \tau_n\,\left|\frac{w^n_r-w^{n-1}_r}
{\tau_n}\right|_{0,\Omega}^2
+
\sum_{n=1}^N \tau_n\,|\del\,.\,\vq^n_r|_{0,\Omega}^2
+  \max_{n=1,\dots, N} |\del w^n_r|_{0,p,\Omega}
+
\sum_{n=1}^N \tau_n\,|\vq^n_r|_{0,r,\Omega}^{2r}
&\leq C.
\label{energy5}
\end{align}
\end{lmm}
\begin{proof}
Choosing $\eta \equiv [\widehat F(w^n_r+b_e(t_n))-\widehat F(w^{n-1}_r+b_e(t_{n-1}))]
-[\widehat F(b_e(t_n))-\widehat F(b_e(t_{n-1}))]$ in
(\ref{Qt1rsc}), and noting (\ref{Qt2rsc}) and (\ref{abcoer}),
yields for $n=2, \ldots, N$
that
\eqlabstart
\begin{align}
&\tau_n\,\left(\frac{(w^n_r+b_e(t_n))-(w^{n-1}_r+b_e(t_{n-1}))}{\tau_n},
\frac{[\widehat F(w^n_r+b_e(t_n))-\widehat F(w^{n-1}_r+b_e(t_{n-1}))]
-[\widehat F(b_e(t_n))-\widehat F(b_e(t_{n-1}))]}{\tau_n}\right) \nonumber \\
& \hspace{1in} = -\left(\del\,.\,\vq^n_r,
[\widehat F(w^n_r+b_e(t_n))-\widehat F(w^{n-1}_r+b_e(t_{n-1}))]-[\widehat F(b_e(t_n))-\widehat F(b_e(t_{n-1}))]
\right)
\nonumber \\
&\hspace{1in}
= -\left(k\,\left[|\vq^n_r|^{r-2}\,\vq^n_r -
|\vq^{n-1}_r|^{r-2}\,\vq^{n-1}_r\right], \vq^n_r \right)
\nonumber \\
& \hspace{1in} \leq - \tfrac{r-1}{r}\,
\left(k,|\vq^n_r|^r-|\vq^{n-1}_r|^r\right),
\label{Fdiff}
\end{align}
and, on noting (\ref{inidata}),
\begin{align}
&\tau_1\,\left(\frac{(w^1_r+b_e(t_1))-(w^0_r+b_e(t_0))}{\tau_1},
\frac{[\widehat F(w^1_r+b_e(t_1))-\widehat F(w^0_r+b_e(t_0))]
-[\widehat F(b_e(t_1))-\widehat F(b_e(t_0))]}{\tau_1}\right) \nonumber \\
& \hspace{1in} = -(\del\,.\,\vq^1_r,
[\widehat F(w^1_r+b_e(t_1))-\widehat F(w^0_r+b_e(t_0))]-[\widehat F(b_e(t_1))-\widehat F(b_e(t_0))])
\nonumber \\
& \hspace{1in}
= -(k,|\vq^1_r|^{r}) - (\vq^1_r,[\widehat M(w_0+b_e(t_0))]^{-1} \, \del w_0)
\nonumber \\
& \hspace{1in}
\leq
-(k,|\vq^1_r|^{r}) + (k, |\vq^1_r|)
\nonumber \\
& \hspace{1in}
\leq  - \tfrac{r-1}{r}\,(k,|\vq^1_r|^{r}) + \tfrac{1}{p} \,|k|_{0,1,\Omega}\,.
\label{Fdiff1}
\end{align}
\eqlabend
Summing (\ref{Fdiff}) and including (\ref{Fdiff1}) 
yields for
$n = 1,\ldots,N$ that
\begin{align}
&\tfrac{r-1}{r}\,
(k\,|\vq^n_r|^r,1)
+ \sum_{\ell=1}^n \tau_\ell\,\left(\frac{(w^\ell_r+b_e(t_\ell))-(w^{\ell-1}_r+b_e(t_{\ell-1}))}{\tau_\ell},
\frac{\widehat F(w^\ell_r+b_e(t_\ell))-\widehat F(w^{\ell-1}_r+b_e(t_{\ell-1}))}{\tau_\ell}\right)
\nonumber \\
&\hspace{1in}\leq C + \sum_{\ell=1}^n \tau_\ell\,
\left(\frac{(w^\ell_r+b_e(t_\ell))-(w^{\ell-1}_r+b_e(t_{\ell-1}))}{\tau_\ell},
\frac{\widehat F(b_e(t_\ell))-\widehat F(b_e(t_{\ell-1}))}{\tau_\ell}\right).
\label{Fdiffsum}
\end{align}
The first two bounds in the desired result (\ref{energy5}) then follow
from (\ref{Fdiffsum}), (\ref{wF}), (\ref{fn}), (\ref{fnsum}) and (\ref{calFn}),
on using a Young's inequality.
The third bound in (\ref{energy5}) then follows from the second bound in (\ref{energy5}),
(\ref{Qt1rsc}) with $\eta= \del\,.\,\vq^n_r$, (\ref{fn}), (\ref{fnsum}) and (\ref{calFn}).

Next, we prove the fourth bound in (\ref{energy5}).
First, we note from (\ref{Qt2rsc}), (A1) and the first bound in (\ref{energy5})
that for $n=1, \ldots, N$
\begin{align}
|(\widehat F(w^n_r + b_e(t_n))-\widehat F(b_e(t_n)), \del\,.\,\vv)|
&\leq |k|_{0,\infty,\Omega}\,|\vq^n_r|_{0,r,\Omega}^{r-1}\,|\vv|_{0,r,\Omega}
\leq \left[ \textstyle \frac{r-1}{r}\,|\vq^n_r|_{0,r,\Omega}^{r}+C\right]
\,|\vv|_{0,r,\Omega}
\nonumber \\
&\leq C\,|\vv|_{0,r,\Omega} \qquad
\forall \ \vv \in \vV^r(\Omega)\,.
\label{gradF}
\end{align}
It follows from (\ref{gradF}),
as $C^\infty_0(\Omega)$ is dense in $L^r(\Omega)$,
that the distributional gradient of $\widehat F(w^n_r+b_e(t_n))$ belongs to the dual of
$[L^r(\Omega)]^2$. 
Hence,
we deduce from (\ref{gradF}) that  
\begin{align}
|\del[\widehat F(w^n_r+b_e(t_n))]|_{0,p,\Omega}
&\leq C,  
\qquad n=1, \ldots, N.
\label{gradFnew}
\end{align}
As $\widehat F$ is globally Lipschitz, recall (\ref{wF}),
we obtain the
fourth bound in (\ref{energy5}).

Finally, it follows from (\ref{Qt2rsc}) with $\vv = \vq^n_r$,
(\ref{jcsc}), (\ref{wF}), (\ref{energy4}) and the third bound in (\ref{energy5}) that
\begin{align}
k_{\rm min}^2 \sum_{n=1}^N \tau_n\,|\vq^n_r|_{0,r,\Omega}^{2r}
\leq \sum_{n=1}^N \tau_n [\,(k,|\vq^n_r|^r)\,]^2
&= \sum_{n=1}^N \tau_n [\,(\hat F(w^n_r+b_e(t_n)-\hat F(b_e(t_n)), \del\,.\,\vq^n_r)\,]^2
\nonumber \\
& \leq (\hat M_0)^{-2}
\sum_{n=1}^N \tau_n \,|w^n_r|_{0,\Omega}^2\,|\del\,.\,\vq^n_r|_{0,\Omega}^2
\leq C\sum_{n=1}^N \tau_n \,|\del\,.\,\vq^n_r|_{0,\Omega}^2 \leq C.
\label{qnr2r}
\end{align}
Hence, the final bound in (\ref{energy5}) holds.
\end{proof}

\subsubsection{Convergence of (Q$^{\tau}_r$) to (Q$_r$)}

In addition to the notation (\ref{wqtau+}),
we introduce for $t \in (t_{n-1},t_n]$, $n = 1, \ldots, N$,
\begin{align}
b_e^{\tau}(t) &
:= \frac{(t-t_{n-1})}{\tau_n}\,b_e(t_n)
+\frac{(t_n-t)}{\tau_n}\,b_e(t_{n-1}), \qquad
b_e^{\tau,+}(t) := b_e(t_n)\,.
 \label{wqtau+be}
\end{align}
We also write $w_r^{\tau(,+)}$ to mean with or without the superscript $+$.
We note from (\ref{wqtau+be}),  (\ref{calFn}), (\ref{fn}) and (\ref{fnsum}) that
\begin{align}
b_e^{\tau},\, b_e^{\tau,+} \rightarrow b_e \quad \mbox{strongly in } L^2(0,T), \qquad
\frac{db_e^{\tau}}{dt} \rightarrow \frac{db_e}{dt} \quad \mbox{strongly in } L^2(0,T)
\quad \mbox{ as } \tau \rightarrow 0.
\label{ftauconv}
\end{align}
We set also $\Omega_T := \Omega \times (0,T)$.

Adopting the notation (\ref{wqtau+}) and (\ref{wqtau+be}),
(Q$^{\tau}_r$), (\ref{Qt1rsc},b),  can be restated as:
Find $w_r^{\tau} \in H^1(0,T;L^2(\Omega))$
and $\vq_r^{\tau,+} \in L^2(0,T;\vV^r(\Omega))$
such that
\eqlabstart
\begin{align}
\int_0^T \left(\frac{\partial w_r^\tau}{\partial t} +
\del\,.\,
\vq^+_r  + \frac{db_e^{\tau}}{dt},\eta\right)\,{\rm d}t &=0
\qquad \forall
\ \eta \in L^2(\Omega_T),
\label{W+-}\\
\int_{0}^T \left[\,
(k\,|\vq^{\tau,+}_r|^{r-2} \vq^{\tau,+}_r,\vv)
-(\widehat F(w^{\tau,+}_r+b_e^{\tau,+})-\widehat F(b^{\tau,+}_e), \del\,.\,\vv)
\right]\,{\rm d}t &=0 \qquad \forall
\ \vv \in L^2(0,T;\vV^r(\Omega));
\label{Ptph+}
\end{align}
\eqlabend
where $w_r^{\tau}(\cdot,0)=w_0(\cdot)$.

In Theorem \ref{conthmscr} below we show the convergence of (Q$^{\tau}_r$),
(\ref{W+-},b), as $\tau \rightarrow 0$ to

\vspace{2mm}

\noindent
{\bf(Q$_r$)}
Find $w_r \in H^1(0,T;L^2(\Omega))$
and $\vq_r \in L^2(0,T;\vV^r(\Omega))$
such that
\eqlabstart
\begin{alignat}{2}
\int_0^T \left(\frac{\partial w_r}{\partial t} + \del \,.\,\vq_r + \frac{d{b}_e}{dt}
,\eta\right) \dt &=0
\qquad &&\forall \ \eta \in
L^2(
\Omega_T)
,
\label{Qaesc}\\
\int_0^T \left[ (k\,|\vq_r|^{r-2} \vq_r, \vv)
- (\widehat F(w_r+b_e)-\widehat F(b_e), \del\,.\,\vv)
\right] \dt &= 0
\qquad &&\forall \ \vv \in
 L^2(0,T;\vV^r(\Omega))\,;
\label{Qbesc}
\end{alignat}
\eqlabend
where $w_r(\cdot,0)=w_0(\cdot)$.

\begin{thrm}\label{conthmscr}
Let the Assumptions (A1) and (A2) hold.
For any fixed  $r \in (1,\frac{4}{3})$
and for all time partitions $\{\tau_n\}_{n=1}^N$,
there exists a subsequence of
$\{w_r^{\tau},\vq^{\tau,+}_r\}_{\tau>0}$ (not indicated),
where $\{w_r^{\tau},\vq^{\tau,+}_r\}$
is the unique solution of (Q$^{\tau}_{r}$), (\ref{W+-},b),
such that as $\tau \rightarrow 0$
\eqlabstart
\begin{alignat}{2}
 w^{\tau(,+)}_r &\rightarrow w_r \qquad &&
\mbox{weak-$\star$ in } L^\infty(0,T;W^{1,p}(\Omega)),
\label{wconscr}
\\
\frac{\partial  w^{\tau}_r}{\partial t} &\rightarrow
\frac{\partial w_r}{\partial t}  \qquad &&
\mbox{weakly in } L^2(\Omega_T),
\label{wtconscr}
\\
w^{\tau(,+)}_r &\rightarrow w_r  \qquad &&
\mbox{strongly in } L^2(\Omega_T),
\label{sconscr}
\\
\vq^{\tau,+}_r &\rightarrow \vq_r \qquad &&
\mbox{weakly in } L^{2r}(0,T;[L^r(\Omega)]^2), \label{qMcon1scr}\\
\del\,.\,\vq^{\tau,+}_r &\rightarrow \del\,.\,\vq_r \qquad &&
\mbox{weakly in } L^2(\Omega_T). \label{delqconscr}
\end{alignat}
\eqlabend
Moreover, $\{w_r,\vq_r\}$ solves (Q$_r$), (\ref{Qaesc},b). 
\end{thrm}
\begin{proof}
The bounds (\ref{energy4}) and (\ref{energy5}) yield
immediately that
\eqlabstart
\begin{align}
&\|w_r^{\tau(,+)}\|_{L^\infty(0,T;W^{1,p}(\Omega))}
+
\left\|\frac{\partial w_r^{\tau}}{\partial t}\right\|_{L^2(\Omega_T)}^2
+ 
\|\vq_r^{\tau,+}\|_{L^{2r}(0,T;L^r(\Omega))}^{2r}
+
\|\del\,.\,\vq_r^{\tau,+}\|_{L^2(\Omega_T)}^2
\leq C,
\label{+boundsscr}
\\
&\|w^{\tau}_r-w^{\tau,+}_r\|_{L^2(\Omega_T)} \leq \tau^2
\left\|\frac{\partial w_r^{\tau}}{\partial t}\right\|_{L^2(\Omega_T)}^2
\leq C\,\tau^2.
\label{timediffsc}
\end{align}
\eqlabend
The subsequence convergence results (\ref{wconscr}--e) follow immediately
from (\ref{+boundsscr},b).
The strong convergence result (\ref{sconscr}) follows from
(\ref{wconscr},b), the compactness result (\ref{compact1}) and (\ref{timediffsc}).
As $w_r^\tau(\cdot,0)=w_0(\cdot)$, it follows from the above that
$w_r(\cdot,0)=w_0(\cdot)$.

It follows immediately from passing to the limit
$\tau \rightarrow 0$ in (\ref{W+-}) for the subsequence,
on noting (\ref{wtconscr},e) and (\ref{ftauconv}),
that $\{w_r,\vq_r\}$ satisfy (\ref{Qaesc}). 

Given any $\vz \in
L^2(0,T;[\vV^r(\Omega)]^2)$, we choose
$\vv \equiv \vq^{\tau,+}_r-\vz$ in (\ref{Ptph+})
to yield, on noting (\ref{aederiv}), that
\begin{align}
\int_{0}^T
(\widehat F(w^{\tau,+}_r+b^{\tau,+}_e)- \widehat F(b^{\tau,+}_e),
\del\,.\,(\vq_r^{\tau,+} - \vz)) \,{\rm d}t
&= \int_{0}^T ( k\,|\vq^{\tau,+}_r|^{r-2} \vq^{\tau,+}_r,\vq^{\tau,+}_r-\vz)
\,{\rm d}t
\nonumber \\
&\geq
\int_{0}^T ( k\,|\vz|^{r-2} \vz,\vq^{\tau,+}_r-\vz)
\,{\rm d}t.
\label{ineq1scr}
\end{align}
Passing to the limit $\tau \rightarrow 0$ in (\ref{ineq1scr}) for the subsequence
yields, on noting (\ref{sconscr}--e), (\ref{wF}) and (\ref{ftauconv}), that
\begin{align}
\int_{0}^T
(\widehat F(w_r+b_e)- \widehat F(b_e),
\del\,.\,(\vq_r - \vz)) \,{\rm d}t
&\geq
\int_{0}^T ( k\,|\vz|^{r-2} \vz,\vq_r-\vz)
\,{\rm d}t.
\label{ineq2scr}
\end{align}
For any fixed $\vv \in \vV^r(\Omega)$, choosing $\vz =\vq_r \pm \alpha \,\vv$ with
$\alpha \in {\mathbb R}_{>0}$ in (\ref{ineq2scr}),
and letting $\alpha \rightarrow 0$
yields the desired result
(\ref{Qbesc}). 
Hence $\{w_r,\vq_r\}$ solves (Q$_r$), (\ref{Qaesc},b).
\end{proof}

\subsubsection{Convergence of (Q$_r$) to (Q)}

We need an extra assumption.

\vspace{2mm}

\noindent
{\bf (A5)} $k \in C(\overline{\Omega})$.

\vspace{2mm}

Then the weak mixed formulation of the cylindrical superconductor problem is:

\noindent
{\bf(Q)}
Find $w \in H^1(0,T;L^2(\Omega))$
and $\vq \in L^2(0,T;\vV^{\cal M}(\Omega))$
such that
\eqlabstart
\begin{alignat}{2}
\int_0^T \left(\frac{\partial w}{\partial t} + \del \,.\,\vq + \frac{d{b}_e}{dt}
,\eta\right) \dt &=0
\qquad &&\forall \ \eta \in
L^2(0,T;L^2(\Omega)),
\label{Qasc}\\
\int_0^T \left[ \langle |\vv|
-|\vq|, k \rangle_{C(\overline{\Omega})}
- (\del\,.\,(\vv-\vq), \widehat F(w+b_e)-\widehat F(b_e)\,)
\right] \dt &\ge 0
\qquad &&\forall \ \vv \in
 L^2(0,T;\vV^{\cal M}(\Omega))\,;
\label{Qbsc}
\end{alignat}
\eqlabend
where $w(\cdot,0)=w_0(\cdot)$.

Recalling (\ref{hatK}) and (\ref{jcsc}), it follows that
\begin{align}
\widehat K(\psi):= \{ \eta \in W^{1,\infty}_0(\Omega) : |\del \eta| \leq k\,
\widehat M(\psi) \;
\mbox{ a.e. in }
\Omega \}.
\label{Ksc}
\end{align}
Associated with the mixed formulation (Q) is the primal variational inequality:

\noindent
{\bf (P)} Find $w \in L^{\infty}(0,T;\widehat K(w+b_e))\cap H^1(0,T;L^2(\Omega))$ such that
\begin{align}
\int_0^T \left(\frac{\partial w}{\partial t} + \frac{d{b}_e}{dt}
,\eta-w\right) \, \dt  &\geq 0
\qquad \forall \ \eta \in
L^2(0,T;\widehat K(w+b_e)),
\label{Psc}
\end{align}
where
$w(\cdot,0)=w_0(\cdot)$.

\begin{thrm}\label{conthmsc}
Let the Assumptions (A1), (A2), (A3) and (A5) hold.
Then there exists a subsequence of
$\{w_r,\vq_r\}_{r\in(1,\frac{4}{3})}$ (not indicated),
where $\{w_r,\vq_r\}$
solves (Q$_{r}$), (\ref{Qaesc},b),
such that as $r \rightarrow 1$
\eqlabstart
\begin{alignat}{2}
 w_r &\rightarrow w \qquad &&
\mbox{weak-$\star$ in } L^\infty(0,T;W^{1,4}(\Omega)),
\label{wconsc}
\\
\frac{\partial  w_r}{\partial t} &\rightarrow
\frac{\partial w}{\partial t}  \qquad &&
\mbox{weakly in } L^2(\Omega_T),
\label{wtconsc}
\\
w_r &\rightarrow w  \qquad &&
\mbox{strongly in } C([0,T];C(\overline{\Omega})),
\label{sconsc}
\\
\vq_r &\rightarrow \vq \qquad &&
\mbox{weakly in } L^2(0,T;[{\cal M}(\overline{\Omega})]^2), \label{qMcon1sc}\\
\del\,.\,\vq_r &\rightarrow \del\,.\,\vq \qquad &&
\mbox{weakly in } L^2(\Omega_T). \label{delqconsc}
\end{alignat}
\eqlabend
Moreover, $\{w,\vq\}$ solves (Q), (\ref{Qasc},b). 
\end{thrm}
\begin{proof}
On noting that $r \in (1,\frac{4}{3}) \Rightarrow p>4$,
the results (\ref{+boundsscr}), (\ref{wconscr}--e) and (\ref{aederiv}) yield
immediately that
\begin{align}
&\|w_r\|_{L^\infty(0,T;W^{1,4}(\Omega))}
+
\left\|\frac{\partial w_r}{\partial t}\right\|_{L^2(\Omega_T)}^2
+
\|\vq_r\|_{L^2(0,T;L^1(\Omega))}
+
\|\del\,.\,\vq_r\|_{L^2(\Omega_T)}^2
\leq C(T).
\label{+boundssc}
\end{align}
The subsequence convergence results (\ref{wconsc},b,d,e) follow immediately
from (\ref{+boundssc}).
The strong convergence result (\ref{sconsc}) follows from
(\ref{wconsc},b) and the compactness result (\ref{compact1}).
As $w_r(\cdot,0)=w_0(\cdot)$, it follows from the above that
$w(\cdot,0)=w_0(\cdot)$.
It follows immediately from passing to the limit
$r \rightarrow 1$ in (\ref{Qaesc}) for the subsequence,
on noting (\ref{wtconsc},e),
that $\{w,\vq\}$ satisfy (\ref{Qasc}).

Given any $\vz \in
L^2(0,T;[C^{\infty}(\overline{\Omega})]^2)$, we choose
$\vv \equiv \vq_r-\vz$ in (\ref{Qbesc})
to yield, on noting (\ref{aederiv}), that
\begin{align}
\int_{0}^T
(\widehat{F}(w_r+b_e)-\widehat{F}(b_e),\del\,.\,(\vq_r - \vz)) \,{\rm d}t
&= \int_{0}^T ( k\,|\vq_r|^{r-2} \vq_r,\vq_r-\vz)
\,{\rm d}t
\geq
r^{-1}
\int_{0}^T (k,
|\vq_r|^r - |\vz|^r)  \,{\rm d}t.
\label{ineq1}
\end{align}
It follows immediately from
(\ref{sconscr},e) and (\ref{wF}) that
for any $\vz \in L^2(0,T;[C^{\infty}(\overline{\Omega})]^2)$
\begin{align}
&\int_{0}^T
(\widehat F(w_r+b_e)-\widehat F(b_e),\del\,.\,(\vq_r^+ - \vz)) \,{\rm d}t
\rightarrow
\int_{0}^T
(\widehat F(w+b_e)-\widehat F(b_e),\del\,.\,(\vq - \vz)) \,{\rm d}t
\qquad \mbox{as} \quad r \rightarrow 1.
\label{ineq2}
\end{align}
Next, we note that for any $\vz \in L^2(0,T;[C^{\infty}(\overline{\Omega})]^2)$
\begin{align}
r^{-1}\,\int_{0}^{T} (k\,|\vz|^r,1) \,{\rm d}t
\rightarrow \int_{0}^T (k\,|\vz|,1)\,{\rm d}t
\qquad \mbox{as} \quad r \rightarrow 1.
\label{ineq3}
\end{align}
Finally, it follows from (\ref{qMcon1sc}),
and similarly to (\ref{Mweak2}), that
\begin{align}
\liminf_{r \rightarrow 1} r^{-1} \int_{0}^{T} (k\,|\vq_r|^r ,1) \,{\rm d}t
\geq
\liminf_{r \rightarrow 1} \int_{0}^{T} (k\,|\vq_r|,1) \,{\rm d}t
\geq
\int_0^T
\langle |\vq|,k \rangle_{C(\overline{\Omega})} \,{\rm d}t.
\label{ineq4}
\end{align}
Combining (\ref{ineq1})--(\ref{ineq4}), it follows  that
$\{w,\vq\}$ satisfies (\ref{Qbsc}) for any $\vv \in L^2(0,T;[C^\infty
(\overline{\Omega})]^2)$.
The desired result, $\{w,\vq\}$ satisfies (\ref{Qbsc}) 
for any $\vv \in L^2(0,T;\vV^{\cal M}(\Omega))$,
and hence $\{w,\vq\}$ solves (Q), (\ref{Qasc},b), 
then follows from the density results (1.22b,c)
in Barrett and Prigozhin \cite{BP5}
with ``$=$''
replaced by ``$\leq$'' in the latter.
\end{proof}

\begin{thrm} \label{eqlem}
Let the assumptions of Theorem \ref{conthmsc} hold.
We then have that any solution $\{w,\vq\}$ of (Q), (\ref{Qasc},b), 
satisfies
\begin{align}
&\int_0^T
\left[
\langle |\vq| ,k\,\widehat M(w+b_e)\rangle_{C(\overline{\Omega})}-(w,\del\,.\,\vq)
\right] \,{\rm d}t =0.
\label{m1a}
\end{align}
Moreover, $w$ solves the quasi-variational inequality (P), (\ref{Psc}).
\end{thrm}
\begin{proof}
See the proof of Theorem 3.3 in Barrett and Prigozhin \cite{BP5}.
However, we note that one can establish (\ref{m1a}) by only requiring
the density results (1.22b,c) in Barrett and Prigozhin \cite{BP5}, with ``$=$''
replaced by ``$\leq$'' in the latter, as opposed to (1.22a--c) there.
To see this, we note that it follows immediately from (\ref{Qbesc}), (\ref{+boundssc})
and (\ref{wF}) that
\begin{align}
\int_0^T
\left[
(k\,\widehat M(w_r+b_e)\,|\vq_r|^{r-2}\vq_r , \vv) - (w_r,\del\,.\,\vv)
\right] \,{\rm d}t =0 \qquad \forall \ \vv \in L^2(0,T;\vV^r(\Omega)).
\label{m1ar}
\end{align}
For any fixed $\vz \in L^2(0,T;[C(\overline{\Omega})]^2)$, we choose $\vv=\vq_r-\vz$
in (\ref{m1ar}) and deduce from (\ref{sconsc}--e),
similarly to (\ref{ineq1})--(\ref{ineq4}),
on passing to the limit $r \rightarrow 1$ that
\begin{align}
\int_0^T
\left[
\langle |\vz| - |\vq|, k\,\widehat M(w+b_e) \rangle_{C(\overline{\Omega})} - (\del\,.\,(\vz-\vq),w)
\right] \,{\rm d}t \geq 0 \qquad \forall \ \vz \in L^2(0,T;[C(\overline{\Omega})]^2).
\label{m2ar}
\end{align}
Applying the stated density results from \cite{BP5}, we obtain
(\ref{m2ar}) holds for all $\vz \in L^2(0,T;\vV^{\cal M}(\Omega))$.
Then choosing $\vz=\vzero$ and $\vz=2\,\vq$ in (\ref{m2ar}) yields the desired
result (\ref{m1a}).
\end{proof}

\subsection{Convergence of (Q$^{\tau}_r$) to (Q) in case (iii)}

It follows from (\ref{Qt1r}), (\ref{energy4}), (\ref{calA}), (\ref{bilforma}) and (\ref{calFn})
in the thin film superconductor case that for $n=1,\ldots,N$
\begin{align}
\tau_n\, | (\vq^n_r, \del \eta) | \leq C\,\|\eta\|_{\Hhalf}
\qquad \forall \ \eta \in C^\infty_0(\Omega).
\label{divqnrbdtf}
\end{align}
Hence, for a fixed time partition $\{\tau_n\}_{n=1}^N$,
the distributional divergence of $\vq^n_r$ belongs $[\Hhalf]^\star$, $n=1,\ldots,N$.
On recalling (\ref{vZs}),
(Q$^\tau_r$), (\ref{Qt1r},b), can then be reformulated for a given $r \in (1, \frac{4}{3})$ as:
\vspace{2mm}

\noindent
{\bf (Q$^{\tau}_r$)}
For $n = 1,\ldots, N$, find $w^n_r \in W^{1,p}_0(\Omega)$ and
$\vq^n_r \in \vZ^r(\Omega)$ such that
\eqlabstart
\begin{alignat}{2}
a\left(\frac{w^n_r-w^{n-1}_r}{\tau_n},\eta\right) + \langle \del\,.\,\vq^n_r,\eta
\rangle_{\Hhalf}
 &= -\left(\frac{b_e(t_n)-b_e(t_{n-1})}{\tau_n}
,\eta\right)
\quad \;&&\forall \ \eta \in \Hhalf,\label{Qt1rtf}\\
(k\,|\vq^n_r|^{r-2}\vq^n_r,\vv) - \langle \del\,.\,\vv, w^n_r \rangle_{\Hhalf}  &=0
\quad \;
&&\forall \ \vv\in \vZ^r(\Omega)\,;
\label{Qt2rtf}
\end{alignat}
\eqlabend
where $w^0_r = w_0$.

We have the following stability result.

\begin{lmm}\label{conthmtf}
Let the Assumptions (A1) and(A2) hold.
For any fixed $r \in (1,\tfrac{4}{3})$ and time partition $\{\tau_n\}_{n=1}^N$,
the unique solution $\{w^n_r,\vq^n_r\}_{n=1}^N$ of (Q$^{\tau}_r$), (\ref{Qt1rtf},b), in addition to satisfying
(\ref{energy4}) with $\|\cdot\|_{\cal A} \equiv \|\cdot\|_{\Hhalf}$ satisfies
\begin{align}
\frac{r-1}{r}\,\max_{n=1,\dots, N}
(k,|\vq^n_n|^r)
+ \sum_{n=1}^N \tau_n\,\left\|\frac{w^n_r-w^{n-1}_r}
{\tau_n}\right\|_{\Hhalf}^2
+
\sum_{n=1}^N \tau_n\,\|\del\,.\,\vq^n_r\|_{[\Hhalf]^\star}^2
+ \sum_{n=1}^N \tau_n\,
|\vq^n_r|_{0,r,\Omega}^{2r}
&\leq C.
\label{energy6}
\end{align}
\end{lmm}
\begin{proof}
Similarly to (\ref{Fdiff},b),
choosing $\eta \equiv w^n_r -w^{n-1}_r$ in
(\ref{Qt1rtf}), and noting (\ref{Qt2rtf}) and (\ref{abcoer}),
yields for $n=2, \ldots, N$
that
\eqlabstart
\begin{align}
&\tau_n\,a\left(\frac{w^n_r-w^{n-1}_r}{\tau_n},
\frac{w^n_r-w^{n-1}_r}{\tau_n}\right)
+\tau_n\,\left(\frac{b_e(t_n)-b_e(t_{n-1})}{\tau_n},
\frac{w^n_r-w^{n-1}_r}{\tau_n}\right)
\nonumber \\
& \hspace{3in}
= -\langle\del\,.\,\vq^n_r,
w^n_r-w^{n-1}_r \rangle_{\Hhalf}
\nonumber \\
&\hspace{3in}
= -\left(k\,\left[|\vq^n_r|^{r-2}\,\vq^n_r -
|\vq^{n-1}_r|^{r-2}\,\vq^{n-1}_r\right], \vq^n_r \right)
\nonumber \\
& \hspace{3in}
\leq - \tfrac{r-1}{r}\,
\left(k,|\vq^n_r|^r-|\vq^{n-1}_r|^r\right),
\label{Fdifftf}
\end{align}
and, on noting (\ref{inidata}),
\begin{align}
\tau_1\,a\left(\frac{w^1_r-w^0_r}{\tau_1},
\frac{w^1_r-w^0_r}{\tau_1}\right) + \tau_1\,\left(\frac{b_e(t_1)-b_e(t_{0})}{\tau_1},
\frac{w^1_r-w^{0}_r}{\tau_1}\right) 
& 
= -\langle \del\,.\,\vq^1_r,
w^1_r-w^0_r\rangle_{\Hhalf}
\nonumber \\
& 
= -(k,|\vq^1_r|^{r}) - (\vq^1_r,\del w_0)
\nonumber \\
& 
\leq
-(k,|\vq^1_r|^{r}) + (k, |\vq^1_r|)
\nonumber \\
& 
\leq  - \tfrac{r-1}{r}\,(k,|\vq^1_r|^{r}) + \tfrac{1}{p} \,|k|_{0,1,\Omega}\,.
\label{Fdiff1tf}
\end{align}
\eqlabend
Summing (\ref{Fdifftf}) and including (\ref{Fdiff1tf})
yields for
$n = 1,\ldots,N$ that
\begin{align}
&\tfrac{r-1}{r}\,
(k\,|\vq^n_r|^r,1)
+ \sum_{\ell=1}^n \tau_\ell\,a\left(\frac{w^\ell_r-w^{\ell-1}_r}{\tau_\ell},
\frac{w^\ell_r-w^{\ell-1}_r}{\tau_\ell}\right)
\leq C + \sum_{\ell=1}^n \tau_\ell\,
\left(\frac{b_e(t_\ell)-b_e(t_{\ell-1})}{\tau_\ell},
\frac{w^\ell_r-w^{\ell-1}_r}{\tau_\ell}\right).
\label{Fdiffsumtf}
\end{align}
The first two bounds in the desired result (\ref{energy6}) then follow
from (\ref{Fdiffsumtf}), (\ref{fn}), (\ref{fnsum}) and (\ref{calFn}),
on using a Young's inequality.
The third bound in (\ref{energy6}) then follows from the second bound in (\ref{energy6})
and (\ref{Qt1rtf}).

Finally, similarly to (\ref{qnr2r}), it follows from (\ref{Qt2rtf}) with $\vv = \vq^n_r$,
(\ref{jcsc}), (\ref{wF}), (\ref{energy4}) and the third bound in (\ref{energy6}) that
\begin{align}
k_{\rm min}^2 \sum_{n=1}^N \tau_n\,|\vq^n_r|_{0,r,\Omega}^{2r}
\leq \sum_{n=1}^N \tau_n [\,(k,|\vq^n_r|^r)\,]^2
&= \sum_{n=1}^N \tau_n [\,\langle \del\,.\,\vq^n_r, w^n_r \rangle_{\Hhalf}\,]^2
\nonumber \\
& \leq
\sum_{n=1}^N \tau_n \,\|w^n_r\|_{\Hhalf}^2\,\|\del\,.\,\vq^n_r\|_{[\Hhalf]^\star}^2
\leq C\sum_{n=1}^N \tau_n \,\|\del\,.\,\vq^n_r\|_{[\Hhalf]^\star}^2 \leq C.
\label{qnr2rtf}
\end{align}
Hence, the final bound in (\ref{energy6}) holds.
\end{proof}

\subsubsection{Convergence of (Q$^{\tau}_r$) to (Q$_r$)}

Adopting the notation (\ref{wqtau+}) and (\ref{wqtau+be}),
(Q$^{\tau}_r$), (\ref{Qt1rtf},b), can be rewritten as:
Find $w_r^{\tau}
\in H^{1}(0,T;\Hhalf)$
and
$\vq^{\tau,+}_r \in 
L^2(0,T; \vZ^r(\Omega))$
such that
\eqlabstart
\begin{alignat}{2}
\int_0^T \left[
a(\frac{\partial w_r^{\tau}}{\partial t},\eta)
+ \langle \del\,.\,\vq_r^{\tau,+},\eta \rangle_{\Hhalf} +
(\frac{d b_e^{\tau}}{dt} ,\eta) \right] \,{\rm d}t &=0
\qquad
&&\forall \ \eta \in L^2(0,T;\Hhalf),\label{Q1tau}\\
\int_0^T \left[ (k\,|\vq^{\tau,+}_r|^{r-2}\, \vq^{\tau,+}_r,\vv)
- \langle \del\,.\, \vv, w^{\tau,+}_r \rangle_{\Hhalf} \right] \,{\rm d}t &=0
\qquad &&\forall \ \vv \in
L^2(0,T; \vZ^r(\Omega))
\,;
\label{Q2tau}
\end{alignat}
\eqlabend
where $w_r^\tau(\cdot,0) = w_0(\cdot)$.

In Theorem \ref{conthmtfr} below we show the convergence of
(Q$^\tau_r$), (\ref{Q1tau},b), as $\tau \rightarrow 0$ to

\noindent
{\bf(Q$_{r}$)}
Find $w_r \in H^1(0,T;H^{\frac{1}{2}}_{00}(\Omega))$
and $\vq_r \in
L^2(0,T; \vZ^r(\Omega))$
such that
\eqlabstart
\begin{alignat}{2}
\int_0^T \left[ a(\frac{\partial w_r}{\partial t},\eta) +
\langle  \del\,. \,\vq_r,
\eta \rangle_{\Hhalf}
+ (\frac{d{b}_e}{dt}
,\eta) \right] \dt &=0
\qquad &&\forall \ \eta \in
L^2(0,T;H^{\frac{1}{2}}_{00}(\Omega)),
\label{Qaetf}\\
\int_0^T \left[
(k \,|\vq_{r}|^{r-2}\,\vq_r,\vv) -
\langle \del\,.\,\vv, w_r \rangle_{\Hhalf} \right]
\dt &= 0 \qquad &&\forall \ \vv \in
L^2(0,T; \vZ^r(\Omega))
\,; \label{Qbetf}
\end{alignat}
\eqlabend
where $w_r(\cdot,0)=w_0(\cdot)$.

Associated with (Q$_r$) is the corresponding generalised $p$-Laplacian
problem for $p \in (4,\infty)$ :

{\bf (P$_p$)}
Find $w_r \in L^{p}(0,T;W^{1,p}_0(\Omega))\cap H^1(0,T;\Hhalf)$ such that
\begin{align}
&\int_0^T \biggl[ a(\frac{\partial w_r}{\partial t},\eta)
+
\left(
\frac{|\del w_r|^{p-2}}{k^{p-1}}\,
\del w_r,
\del \eta\right)
+ \left( \frac{d{b}_e}{dt}
,\eta\right) \,\biggr] \dt =0
\qquad \forall \ \eta \in
L^p(0,T;W^{1,p}_0(\Omega)),
\label{Pp}
\end{align}
where $w_r(\cdot,0)=w_0(\cdot)$.

\begin{thrm}\label{conthmtfr}
Let the Assumptions (A1) and (A2) hold.
For any fixed $r \in (1,\frac{4}{3})$
the
sequence $\{ w^\tau_r,\vq^{\tau,+}_r\}_{\tau>0}$,
where $\{w^\tau_{r},\vq^{\tau,+}_{r}\}$ is the unique solution of (Q$^{\tau}_{r}$),
is such that as $\tau \rightarrow 0$
\eqlabstart
\begin{alignat}{2}
w^{\tau(,+)}_{r}&\rightarrow w_r  \qquad &&
\mbox{weak-$\star$ in } L^\infty(0,T;\Hhalf),
\label{Swconr}
\\
\frac{\partial w^{\tau}_{r}}{\partial t}&\rightarrow \frac{\partial w_r}{\partial t}  \qquad &&
\mbox{weakly in } L^2(0,T;\Hhalf),
\label{Stwconr}
\\
\vq^{\tau,+}_r &\rightarrow \vq_r \qquad &&
\mbox{weakly in } L^{2r}(0,T;[L^r(\Omega)]^2),
\label{Vwconr}\\
\del\,.\,\vq^{\tau,+}_r &\rightarrow \del\,.\,\vq_{r} \qquad &&
\mbox{weakly in } L^2(0,T;[\Hhalf]^\star);
\label{Vdivwconr}
\end{alignat}
\eqlabend
where $\{w_r,\vq_r\}$ is the unique solution of (Q$_r$), (\ref{Qaetf},b).
In addition, $w_r$ is the unique solution of (P$_p$), (\ref{Pp}).
\end{thrm}
\begin{proof}
It follows immediately from (\ref{energy4}), (\ref{energy6}) and (\ref{wqtau+}) that
\eqlabstart
\begin{align}
&\|w^{\tau(,+)}_r\|_{L^\infty(0,T;\Hhalf)} +
\|\frac{\partial w^{\tau}_r}{\partial t}\|_{L^2(0,T;\Hhalf)}^2
+ \|\vq^{\tau,+}_r\|_{L^{2r}(0,T;L^r(\Omega))}^{2r}
+ \|\del\,.\,\vq^{\tau,+}_r\|_{L^2(0,T;[\Hhalf]^\star)}^2 \leq C,
\label{contrbds} \\
&\|w^{\tau}_r-w^{\tau(,+)}_r\|_{L^2(0,T;\Hhalf)}^2
\leq \tau^2 \,\|\frac{\partial w^{\tau}_r}{\partial t}\|_{L^2(0,T;\Hhalf)}^2
\leq C\,\tau^2.
\label{timediff}
\end{align}
\eqlabend
It follows immediately from (\ref{contrbds},b) that the 
results (\ref{Swconr}--d)
hold for a subsequence of $\{w^\tau_r,\vq^{\tau,+}_r\}_{\tau>0}$.
We then pass to the limit $\tau \rightarrow 0$ in (\ref{Q1tau}) for the above
subsequence
and obtain (\ref{Qaetf})
for any fixed $\eta \in L^2(0,T;\Hhalf)$,
on noting (\ref{Stwconr},d) and (\ref{ftauconv}).

For any fixed $\veta \in L^2(0,T;\vZ^r(\Omega))$, 
we choose $\vv = \vq^{\tau,+}_r -\vz$ in (\ref{Q2tau}). On noting (\ref{Q1tau}), (\ref{aederiv})
and (\ref{timediff}), we deduce that
\begin{align}
- \int^T_0 \langle \del\,.\,\vz, w^{\tau,+}_r \rangle_{\Hhalf} \,\dt
& = - \int^T_0 \langle \del\,.\,\vq^{\tau,+}_r, w^{\tau,+}_r \rangle_{\Hhalf} \,\dt
+ \int_0^T (k\,|\vq^{\tau,+}_r|^{r-2}\,\vq^{\tau,+}_r,\vq^{\tau,+}_r -\vz) \,\dt
\nonumber \\
&  \geq 
\int^T_0 \left[\, a(\frac{\partial w^{\tau}_r}{\partial t}, w^{\tau,+}_r )
+ (\frac{d b_e^\tau}{\dt}, w^{\tau,+}_r ) \,\right]
\,\dt
+ \int_0^T (k\,|\vz|^{r-2}\,\vz,\vq^{\tau,+}_r -\vz) \,\dt
\nonumber \\
& \geq \ts \frac{1}{2} \ds
\left[ \|w^{\tau}_r(\cdot,T)\|_a^2  - \|w_0\|_a^2 \right) +
\int^T_0 a(\frac{\partial w^{\tau}_r}{\partial t}, w^{\tau,+}_r - w^{\tau}_r)
\, \dt \nonumber \\
& \qquad
+
\int_0^T
( \frac{d b_e^\tau}{dt}, w^{\tau,+}_r )
\,\dt
+ \int_0^T (k\,|\vz|^{r-2}\,\vz,\vq^{\tau,+}_r -\vz) \,\dt
\nonumber \\
& \geq \ts \frac{1}{2} \ds
\left( \|w^{\tau}_r(\cdot,T)\|_a^2  - \|w_0\|_a^2 \right)
+ 
\int_0^T
(\frac{d b_e^\tau}{dt}, w^{\tau,+}_r )
\,\dt
+ \int_0^T (k\,|\vz|^{r-2}\,\vz,\vq^{\tau,+}_r -\vz) \,\dt
- C\,\tau.
\label{conv1}
\end{align}
It follows from (\ref{contrbds}) that
\begin{align}
\|w^{\tau}_r\|_{C([0,T];\Hhalf)} \leq C.
\label{contSig}
\end{align}
Hence we deduce from (\ref{contSig}), on extraction of a possible further subsequence, that
\begin{align}
\liminf_{\tau \rightarrow 0} \|w^{\tau}_r(\cdot,T)\|_a^2 \geq \|w_r(\cdot,T)\|_a^2.
\label{liminf}
\end{align}
On noting (\ref{Swconr},c),
(\ref{liminf}) and (\ref{ftauconv}), we can pass to the
limit $\tau \rightarrow 0$ in (\ref{conv1}) for the above subsequence to obtain
\begin{align}
- \int_0^T \langle \del\,.\,\vz, w_r\rangle_{\Hhalf}  \,\dt &\geq
\ts \frac{1}{2} \ds \left(
\|w_r(\cdot,T)\|_a^2 - \|w_0\|_a^2 \right)
+ \int_0^T (\frac{d b_e}{dt},w_r) \,\dt
+ \int_0^T (k \,|\vz|^{r-2} \vz, \vq_r-\vz ) \,\dt
\nonumber \\
& \hspace{3in}
\qquad \forall \ \vz \in
L^2(0,T;\vZ^r(\Omega)).
\label{conv2}
\end{align}
It follows from (\ref{conv2}) 
and (\ref{Qaetf})
that
\begin{align}
\int_0^T \langle \del\,.\,(\vq_r-\vz), w_r \rangle_{\Hhalf} \,\dt &\geq
\int_0^T (k \,|\vz|^{r-2} \vz, \vq_r-\vz ) \,\dt
\qquad \forall \ \vz \in
L^2(0,T;\vZ^r(\Omega)).
\label{conv3}
\end{align}
For any fixed $\vv \in 
L^2(0,T;\vZ^r(\Omega))$,
choosing $\vz = \vq_r \pm \alpha\,\vv$
with $\alpha \in {\mathbb R}_{>0}$ in (\ref{conv3}), and letting $\alpha \rightarrow 0$
yields the desired result (\ref{Qbetf}).
Hence $\{w_r,\vq_r\}$ solves (Q$_r$), (\ref{Qaetf},b).

In addition, we obtain from (\ref{Q2tau})
for any fixed $\vv \in C([0,T];[C^\infty(\overline{\Omega})]^2)$,
 on noting the third bound in (\ref{contrbds}), that
\begin{align}
\int_0^T (w^{\tau,+}_r, \del\,.\,\vv)  \,\dt
&= \int_0^T (k\,|\vq^{\tau,+}_r|^{r-2}\,\vq^{\tau,+}_r, \vv) \,\dt
\leq C \int_0^T |\vq^{\tau,+}_r|^{r-1}_{0,r,\Omega}
\,|\vv|_{0,r,\Omega}\,\dt
\leq C(T)\,\| \vv\|_{L^r(0,T;L^r(\Omega))}.
\label{GradSrbd}
\end{align}
Passing to the limit $h,\,\tau \rightarrow 0$ in (\ref{GradSrbd}), on noting
(\ref{Swconr}), yields that
\begin{align}
\int_0^T (w_r, \del\,.\,\vv\,) \,\dt
& \leq C\,\|\vv\|_{L^r(0,T;L^r(\Omega))}
\qquad \forall \ \vv \in C([0,T];[C^\infty(\overline{\Omega})]^2).
\label{Gradsrbd}
\end{align}
Hence it follows that
\begin{align}
w_r \in L^p(0,T;W^{1,p}_0(\Omega))
\qquad\mbox{and}\qquad \|w_r\|_{L^{p}(0,T;W^{1,p}(\Omega))} \leq C.
\label{srLp}
\end{align}

To show the uniqueness of this solution  $\{w_r,\vq_r\}$ of (Q$_r$), (\ref{Qaetf},b), and
that $w_r$ is the unique solution of (P$_p$), (\ref{Pp}), see the proof of Theorem 3.1
in Barrett and Prigozhin \cite{EFILMMATH}.
\end{proof}

\subsubsection{Convergence of (Q$_r$) to (Q)}

On recalling (\ref{vZM}) and assuming (A5), 
the weak mixed formulation of the thin film superconductor problem is:

\noindent
{\bf(Q)}
Find $w \in H^1(0,T;H^{\frac{1}{2}}_{00}(\Omega))$
and $\vq \in 
L^2(0,T;\vZ^{\cal M}(\Omega))$
such that
\eqlabstart
\begin{alignat}{2}
\int_0^T \left[ a(\frac{\partial w}{\partial t},\eta) +
\left \langle  \del \,.\,\vq,
\eta\right\rangle_{\Hhalf}
+ \left(\frac{d{b}_e}{dt}
,\eta\right) \right] \dt &=0
\qquad &&\forall \ \eta \in
L^2(0,T;H^{\frac{1}{2}}_{00}(\Omega)),
\label{Qatf}\\
\int_0^T
\left[\langle |\vv| 
-|\vq|, k \rangle_{C(\overline{\Omega})} 
- \left\langle\del\,.\,\left(\vv-\vq\right),w \right\rangle _{\Hhalf}
\right] \dt &\ge 0
\qquad &&\forall \ \vv \in
L^2(0,T;\vZ^{\cal M}(\Omega))\,;
\label{Qbtf}
\end{alignat}
\eqlabend
where $w(\cdot,0)=w_0(\cdot)$.

Let
\begin{align}
K_k:= \{ \eta \in W^{1,\infty}_0(\Omega) : |\del \eta| \leq k \;
\mbox{ a.e. in }
\Omega \}.
\label{Ktf}
\end{align}
Associated with the mixed formulation (Q) is the primal variational inequality:

\noindent
{\bf (P)} Find $w \in L^{\infty}(0,T;K_k)\cap H^1(0,T;\Hhalf)$ such that
\begin{align}
\int_0^T \left[\, a(\frac{\partial w}{\partial t},\eta-w)
+ \left(\frac{d{b}_e}{dt}
,\eta-w\right) \,\right] \dt  &\geq 0
\qquad \forall \ \eta \in
L^2(0,T;K_k),
\label{Ptf}
\end{align}
where
$w(\cdot,0)=w_0(\cdot)$.

\begin{thrm}\label{conthmtfM}
Let the Assumptions (A1), (A2), (A3) and (A5) hold.
Then there exists a subsequence
of $\{w_r,\vq_r\}_{r\in(1,\frac{4}{3})}$, (not indicated),
where $\{w_{r},\vq_{r}\}$ is the unique solution of (Q$_{r}$),
such that as $r \rightarrow 1$
\eqlabstart
\begin{alignat}{2}
w_{r}&\rightarrow w  \qquad &&
\mbox{weak-$\star$ in } L^\infty(0,T;\Hhalf),
\mbox{ weakly in } L^2(0,T;H^1_0(\Omega)),
\label{swcon}
\\
\frac{\partial w_{r}}{\partial t}&\rightarrow \frac{\partial  w}{\partial t}  \qquad &&
\mbox{weakly in } L^2(0,T;\Hhalf),
\label{stwcon}
\\
\vq_r &\rightarrow \vq \qquad &&
\mbox{weakly in } L^2(0,T;[{\cal M}(\overline{\Omega})]^2),
\label{vwcon}\\
\del\,.\,\vq_r &\rightarrow \del\,.\,\vq \qquad &&
\mbox{weakly in } L^2(0,T;[\Hhalf]^\star);
\label{vdivwcon}
\end{alignat}
\eqlabend
where $\{w,\vq\}$ solves (Q), (\ref{Qatf},b).
In addition, $w$ is unique; and the possible non-uniqueness in $\vq$
is restricted to the following: If there were two solutions $\vq^i$, $i=1,\,2$,
then
\begin{align}
\del\,.\,(\vq^2-\vq^1) = 0 \quad a.e.\ {\rm in}\ \Omega_T
\qquad \mbox{and} \qquad \int_0^T \langle |\vq^2|,k\rangle_{C(\overline{\Omega})}\,{\rm d}t
= \int_0^T \langle |\vq^1|,k\rangle_{C(\overline{\Omega})}\,{\rm d}t.
\label{vuniq}
\end{align}
Finally, $w$ is the unique solution of (P), (\ref{Ptf}).
\end{thrm}
\begin{proof}
See the proof of Theorem 3.2 in Barrett and Prigozhin \cite{EFILMMATH}.
\end{proof}

\section{Numerical Algorithm and Simulation Results\label{numexpts}}
\setcounter{equation}{0}
Our iterative procedure for solving the $n^{\rm th}$ step of (Q$^{h,\tau}_{r}$),
(\ref{Qthr1B},b), for $W^n_r$ and $\vQ^n_r$ is as follows.

Set $W^{n,0}_{r} = W^{n-1}_{r}\in N^h_0$ and
$\vQ^{n,0}_{r} = \vQ^{n-1}_{r}\in \vS^h$. For $m \geq 1$, given iterates
$W^{n,m-1}_{r}\in N^h_0$ and $\vQ^{n,m-1}_{r} \in \vS^h$, we use
 the following linearized version of (\ref{Qthr2B})
\begin{align}
{\mathfrak M}^{h,n}(P^h W^{n,m-1}_r)
\left[|\vQ^{n,m-1}_{r}|_{\delta}^{r-2} \,\vQ^{n,m}_{r}
- (\,|\vQ^{n,m-1}_{r}|^{r-2}_{\delta}
-|\vQ^{n,m-1}_{r}|^{r-2})\,\vQ^{n,m-1}_{r}
\right] + \del_h W^{n,m}_r = \vzero,
\label{linit}
\end{align}
where $|\vv|_{\delta} := (|\vv|^2 +\delta^2)^{\frac{1}{2}}$ with $\delta^2 \ll 1$, to obtain
\begin{align}
\vQ^{n,m}_{r}=\vQ^{n,m-1}_{r}-\frac{  |\vQ^{n,m-1}_{r}|^{r-2}\vQ^{n,m-1}_{r}+
\left[{{\mathfrak M}^{h,n}(P^h W^{n,m-1}_r)}\right]^{-1}{\del_h W^{n,m}_r}}
{|\vQ^{n,m-1}_{r}|_{\delta}^{r-2}}.\label{Qit}
\end{align}
Substituting (\ref{Qit}) into an iterative version of (\ref{Qthr1B}), yields the following linear
system for $W^{n,m}_r\in N^h_0$
\begin{align}
&{\cal A}_h\left(\frac{W^{n,m}_r-W^{n-1}_r}{\tau_n}, \eta^h\right)
+ \left(
[{\mathfrak M}^{h,n}(P^h W^{n,m-1}_r)]^{-1}\,|\vQ^{n,m-1}_{r}|_{\delta}^{2-r}
\,\del_h W^{n,m}_r, \del_h \eta^h \right)
\nonumber \\
& \qquad \qquad =
({\cal F}^n,\eta^h) +
\left(
|\vQ^{n,m-1}_{r}|_{\delta}^{2-r}
\left[
\,|\vQ^{n,m-1}_{r}|^{r-2}_{\delta}
-|\vQ^{n,m-1}_{r}|^{r-2}\right]\,\vQ^{n,m-1}_{r}
,\del_h \eta^h\right)
\qquad
\forall \ \eta^h\in N^h_0\,.
\label{Qthr2it}
\end{align}
Clearly, the linear system (\ref{Qthr2it}) is well-posed.
Solving it for $W^{n,m}_r \in N^h_0$, we then obtain $\vQ^{n,m}_r\in \vS^h$ from (\ref{Qit}).

Prior to the next iteration, $\vQ^{n,m}_{r}$ is then replaced by 
$\alpha \,\vQ^{n,m}_{r}+(1-\alpha)\,\vQ^{n,m-1}_{r}$, where $\alpha \in (0,1)$ (under-relaxation) 
was sometimes needed for convergence in case (i) and $\alpha>1$ (over-relaxation) 
led to acceleration of convergence in cases (ii) and (iii).
Although we have no convergence proof of this procedure, in practice it worked well. 
In particular, the number of iterations was almost independent of the mesh size and the value of $r$.
We note that similar algorithms have been used in \cite{sandqvi,BP5,SUST,EFILMMATH}, 
but there a linear system of similar size to (\ref{Qthr2it}) 
was solved on each iteration for the dual variable $\vQ^{n,m}_r$
and then $W^{n,m}_r$ was updated explicitly.

Being a solution to the primal quasi-variational inequality (P), 
the primal variable $w$ is rate-independent. It can be shown, similarly to 
\cite[Section 4]{EFILMMATH}, that if the 
direction of the dual variable $\vq$ does not change with time a.e.\ in the incident set 
$|\nabla w|={\mathfrak M}(w)$ and this set  increases   monotonically in time, 
then the primal variable 
$w$ at time $t$ depends soley on $w^0$ and $\int_0^t {\cal F}(\cdot,s) {\rm d}s$.   
These conditions are satisfied in our examples below. However, the dual variable $\vq$ is not 
rate-independent. Hence, our time step strategy for approximating both $w$ and $\vq$ 
at time $T$, on assuming that they are changing gradually with time, 
was to 
choose a large time step
$\tau_1$ followed by a small time step $\tau_2=T-\tau_1$. 
Then $W^2_r(\cdot)$ is regarded as an approximation to $w(\cdot,T)$,
whereas $\vQ^2_r(\cdot)$ can be
regarded as an approximation to either the mean of $\vq(\cdot,t)$ over 
the time interval $(T-\tau_2,T)$ or, 
as we did in this work, to $\vq(\cdot,T-0.5\,\tau_2)$.

As in \cite{sandqvi}, for ease of implementation in case 
(i) we replaced $w^\epsilon_0$ by $w_0$ in (\ref{inidata}) and (\ref{weps0h}).
Throughout, we set $\delta= 10^{-10}$, chose $r=1+10^{-9}$, and adopted the stopping criterion
\begin{align}
&\frac{  | \pi^h_N [ \,|W_r^{n,m} - W_r^{n,m-1}| \,]\, |_{0,1,\Omega}}
{ | \pi^h_N [\, |W_r^{n,m}|\, ]\, |_{0,1,\Omega}}
<10^{-6} \quad \;\mbox{and}\quad \;
\frac{  | \vQ_r^{n,m} - \vQ_r^{n,m-1} |_{0,1,\Omega}}
{ | \vQ_r^{n,m} |_{0,1,\Omega}} <2\cdot 10^{-5}\,.
\label{stopA}
\end{align}
The simulations have been performed in Matlab R2012b (64 bit) on a PC with
an Intel Core i5-2400 3.1 GHz processor and 8Gb RAM. The Matlab PDE Toolbox was used
for the triangulation of $\Omega$, which was quasi-uniform.
Although for the convergence analysis in the previous sections,
we assumed, for ease of exposition, that $\Omega$ was polygonal and that
the bilinear form  $c(\del_h\cdot,\del_h\cdot)$ on $N^h_0 \times N^h_0$
was calculated exactly; in practice
curved domain boundaries were approximated by polygonal ones
and $c(\del_h\cdot,\del_h\cdot)$ was approximated, see the Appendix in \cite{SUST} for details.

To  compare our  nonconforming approximations (Q$^{h,\tau}_r$), (\ref{Qthr1B},b), 
with those in
\cite{sandqvi,BP5,EFILMMATH} based on the Raviart--Thomas element, 
we considered three problems with known analytical
solutions.

Our first example is a sandpile growing upon the initial support surface $w_0=\max(0.4-|\vx|,0)$ 
below the source $f$, which was uniform in its support $|\vx|\le 0.2$ with $\int_{\Omega}f\,\dx=1$. 
Due to the radial symmetry, the analytical solution to the unregularized problem is easily found, 
see\cite{sandqvi}. We approximated the regularized (with $\epsilon=0.01$ in (\ref{Meps})) 
quasi-variational inequality problem in the square $\Omega=(-1,1)^2$,
with the internal friction of sand $k_0=0.4$, 
both by (Q$^{h,\tau}_r$), (\ref{Qthr1B},b), proposed in this work 
and by the Raviart--Thomas approximation in \cite{sandqvi}. 
In both cases two time steps, $\tau_1=0.19$ and $\tau_2=0.01$, 
were made to obtain the approximation at $T=0.2$. 
On recalling (\ref{Ph}), we estimated the relative errors by
$$\delta(w):=\frac{| P^h\, W^2_r-w^*|_{0,1,\Omega}}{|w^*|_{0,1,\Omega}}\qquad 
\mbox{and}\qquad \delta(\vq):=\frac{|\vQ^2_r-\vq^*|_{0,1,\Omega}}{|\vq^*|_{0,1,\Omega}},$$
for two meshes with $h=0.04$ and $h=0.02$.
Here $w^* \in S^h$ and $\vq^* \in \vS^h$
with $w^* |_{\sigma} = w(\vx^\sigma,0.2)$ and $\vq^* |_{\sigma}
= \vq(\vx^{\sigma},0.195)$,
where $\vx^{\sigma}$ is the centroid of triangle $\sigma \in {\cal T}^h$. 
For the Raviart--Thomas approximation, \cite{sandqvi}, the best convergence was achieved with 
$\alpha=0.7$ (under-relaxation), whilst for (Q$^{h,\tau}_r$) no relaxation was needed with 
the fastest convergence being for $\alpha=1$. 
Although more iterations at each time step were needed, the latter method produced 
a more accurate approximation, 
see Table \ref{Table1}, and was much simpler to realize.
\begin{table}[h!]%
\caption{Growing sandpile. Approximation by the Raviart--Thomas (RT) and \newline 
the nonconforming linear (NC) element}
\begin{tabular}{ccccc}
\hline%
\multirow{2}*{$h$}  & finite &$\delta(w)$ & $\delta(\vq)$ & CPU time \\
     & element&        \%        &    \%     & (min)\\
\hline
\multirow{2}*{0.04} & RT     & 0.38         &  5.0            & 1.1\\
  & NC     & 0.26         &  4.3            & 1.1\\
\hline%
\multirow{2}*{0.02} & RT     & 0.14         &  2.5            & 5.7\\
  & NC     & 0.08         &  2.3            & 6.8\\
\hline%
\end{tabular}\label{Table1}
\end{table}

As our second example, let us consider a cylindrical superconductor 
and assume the Kim critical state model with
$j_c(b)=(1+|b|/B_0)^{-1}$ with $B_0 \in {\mathbb R}_{>0}$, where we recall (\ref{jcsc}). 
Let $w_0=0$ and the external field grow monotonically, $db_e/dt\ge 0$, with $b_e(0)=0$. 
Then the magnetic field in the superconductor, $b(\vx,t)=w(\vx,t)+b_e(t)$, 
can be found analytically, see \cite{BP5}.  At any point in time, this field is a function of the 
distance $d(\vx):=\mbox{dist}(\vx,\partial\Omega)$ to the domain boundary:
$b(\vx,t)=[u(d(\vx),t)]_+$, 
where $s_+:=\max\{s,0\}$ and $u(d,t)$ satisfies $\partial u/ \partial d=-j_c(u)$ with 
$u(0,t)=b_e(t)$. Solving this equation, we obtain that
\begin{equation}
b(\vx,t)=-B_0+\sqrt{B_0^2+2B_0\,[d_0(t)-d(\vx)]_+}\,,\label{bcyl}\end{equation}
where $d_0(t)=b_e(t)(1+0.5\,b_e(t)/B_0)$ is the depth of the field penetration zone at time $t$. 
The current density $\vj$ is critical in the penetration zone, $|\vj(\vx,t)|=j_c(b(\vx,t))$ 
for $d(\vx)\le d_0(t)$, and zero outside of it. 
As $\vj = \del \times w = \del \times b$,
the current streamlines are the level contours of $b$.
It is more difficult to find the electric field for a general domain $\Omega$ but, 
if $\Omega$ is a rectangle, the analytical solution for $\ve$ can be found 
in Brandt \cite{Brandt} for the Bean model, and can be easily extended 
to the Kim model with a field dependent critical current density.

\begin{figure}[!ht]
\begin{floatrow}
\ffigbox{%
  \includegraphics[width=3.3cm]{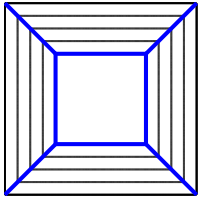}
}{%
  \caption{Current streamlines (thin black) and current density discontinuity lines (thick blue).}%
}
\capbtabbox{%
\begin{tabular}{ccccc}
\hline%
\multirow{2}*{$h$}  & finite &$\delta(w)$ & $\delta(\vq)$ & CPU time \\
     & element&        \%        &    \%     & (min)\\
\hline
\multirow{2}*{0.02} & RT     & 0.25         &  3.5            & 6.2\\
  & NC     & 0.15         &  3.5            & 0.5\\
\hline%
\multirow{2}*{0.01} & RT     & 0.07         &  2.0            & 117\\
  & NC     & 0.05         &  1.9            & 2.8\\
\hline%
\end{tabular}\label{Table2}
}{%
\caption{Cylindrical superconductor.} 
}
\end{floatrow}
\end{figure}

If $\Omega=(0,1)^2$, the field penetration zone consists of four regions of unidirectional 
current density, see Figure 1. Current discontinuity lines separate these regions from each other
and the central zero current region. 
Noting that the direction of the electric field should coincide with that of the current density 
and, as $\nabla \times \ve = - db_e/dt$, the tangential component of $\ve$ 
must be continuous along these discontinuity lines, it follows  
that the electric field should vanish on these lines. 
Let $R_1(t):= \{\vx \in \Omega: x_1 \in (0,1),\; x_2\in (0,s(x_1,t)\}$, where
$x_2=s(x_1,t):=\min(x_1,d_0(t),1-x_1)$ for $x_1\in [0,1]$ is part of the discontinuity lines. 
We have that $\ve=[e_{1}(\vx,t),0]^\top$ in $R_1(t)$ and $e_{1}(x_1,s(x_1,t),t)=0$ for 
$x_1\in [0,1]$. 
Faraday's law, which in $R_1(t)$ reduces to
$\partial e_1/\partial x_2=\partial b/\partial t$, and (\ref{bcyl}) yield that 
$$e_{1}(\vx,t)=\left[\sqrt{B_0^2+2B_0\,(d_0(t)-s(x_1,t))}-\sqrt{B_0^2+2B_0\,(d_0(t)-x_2)}\right]
\,\frac{d}{dt}d_0(t)
\qquad \mbox{in } R_1(t).$$
Similarly, one can find the electric field in the three other regions of the penetration zone. 
Solving the problem numerically, we chose $B_0=0.05$, $b_e(t)=t$ 
and used two time steps, $\tau_1=0.09$ and $\tau_2=0.01$ to find the numerical solution at 
$T=0.1$; see Table \ref{Table2} for a comparison 
of (Q$^{h,\tau}_r$), (\ref{Qthr1B},b),
to the method in \cite{BP5} based on the 
Raviart--Thomas element.
For both methods, over-relaxation with $\alpha=1.8$ led to the fastest convergence. 
The finite element scheme in \cite{BP5} was based on the modified formulation (\ref{Qrv2a}) of (\ref{Qrv2})
which, probably, 
was less efficiently realized in our program. 
This could be the reason for vast difference 
in computation times of the two methods in this case, even though less iterations were 
needed for the method in \cite{BP5}. 
However, the programming of this scheme is more involved, and 
the computed primal variable is less accurate.

Our last example is the magnetization of a thin superconducting disc. For the Bean model, 
$j_c\equiv 1$, the sheet current density and the magnetic field are known, see  
\cite{MikhKuz,ClemSanchez}. Using this analytical solution, the electric field can also be 
calculated, see \cite{SUST}. 
The primal variable, the magnetization function $w$, in thin film magnetization problems is an 
auxiliary variable. Of main interest in such problems are the sheet current density 
$\vj=\del \times w$ and the electric field $\ve$. 
In addition, the magnetic field can be determined from $\vj$ 
by means of the Biot--Savart law, (\ref{BS}). 
To compare the nonconforming approximation (Q$^{h,\tau}_r$), (\ref{Qthr1B},b),
with
the Raviart--Thomas approximation   
in \cite{SUST,EFILMMATH} we present the numerical errors for the two 
main variables, $\delta(\vj)$ and $\delta(\ve)=\delta(\vq)$ in Table \ref{Table3},
where $\delta(\vj)$ is defined similarly to $\delta(\vq)$. 
Since the bilinear form $c(\del_h \cdot,\del_h \cdot)$ on $N^h_0 \times N^h_0$ 
leads to a dense matrix, the numerical solution of (\ref{Qthr2it})
is both memory and time consuming for fine meshes. 
We note that the computation times in Table 3 do not include the time 
for 
assembling the entries of $c(\del_h \cdot,\del_h \cdot)$ on $N^h_0 \times N^h_0$.
Here we recall that these entries were approximated, see the Appendix in \cite{SUST} for details. 
In this example we chose $\Omega$ to be the unit disc, 
$b_e(t)=t$, and found the numerical solution at $T=0.65$ using two time steps, $\tau_1=0.6$ and 
$\tau_2=0.05$. Over-relaxation with $\alpha=1.8$ was employed in both iterative procedures.

\begin{table}[!ht]%
\caption{Thin film magnetization.}
\begin{tabular}{ccccc}
\hline%
\multirow{2}*{$h$}  & finite &$\delta(\vj)$ & $\delta(\vq)$ & CPU time \\
                    & element&        \%        &    \%     & (min)\\
\hline
\multirow{2}*{0.06} & RT     & 0.89         &  3.3            & 4.1\\
                    & NC     & 0.15         &  0.31            & 2.4\\
\hline%
\multirow{2}*{0.03} & RT     & 0.46         &  1.3            & 125\\
                    & NC     & 0.06         &  0.24            & 164\\
\hline%
\end{tabular}\label{Table3}
\end{table}

For the approximation in \cite{SUST,EFILMMATH}, employing the lowest order Raviart--Thomas element 
for $\vq$
and the 
continuous piecewise linear element for $w$,
the approximate current density was calculated directly as $\vJ^n_r=\del \times W^n_r\in \vS^h$.
The same approach was used here for the nonconforming approximation on each element 
$\sigma \in {\cal T}^h$.
However, we note that such a simple procedure may lead to an inaccurate approximation of $\vj$
in thin film problems
involving transport currents, which lead to  non-homogenous time-dependent boundary data for $w_r$
and singular time-dependent forcing data ${\cal F}$ in (\ref{Qrv1},b).
Problems of this type 
have been approximated using the appropriately modified nonconforming
approximation (Q$^{h,\tau}_r$), (\ref{Qthr1B},b), in \cite{Transport}.
There, on recalling (\ref{hatWnr}) and (\ref{wsconhatWnr},b),  
instead of setting $\vJ^n_r=\del \times W^n_r$ on each $\sigma \in {\cal T}^h$,
we set $\vJ^n_r=\del \times  \hat W^n_r\in \vS^h$ and this led to a more accurate approximation 
of $\vj$.  
We note that the cost of the postprocessing step (\ref{hatWnr})
is negligible compared to solving (\ref{Qthr2it}).

\end{document}